\documentclass[11pt,letterpaper]{article}

\usepackage{amsthm}
\usepackage{etoolbox}
\usepackage[margin=1in]{geometry}
\usepackage{amsmath}
\usepackage{graphicx}
\usepackage{booktabs}
\usepackage{amssymb}
\usepackage[most]{tcolorbox}
\usepackage[authoryear,round]{natbib}

\newtcolorbox{graytheorembox}{
    enhanced,
    breakable,
    colback=gray!10,
    colframe=gray!50,
    boxrule=0.5pt,
    arc=2pt,
    left=6pt,
    right=6pt,
    top=6pt,
    bottom=6pt
}

\usepackage{xcolor}
\usepackage{tikz,float,needspace}
\usetikzlibrary{arrows.meta}
\usepackage[
    colorlinks=true,
    citecolor=blue!80!black,
    linkcolor=red,
    urlcolor=teal
]{hyperref}

\usepackage[capitalise]{cleveref}

\usepackage{todonotes}
\usepackage[ruled,vlined]{algorithm2e}

\usepackage{mathtools}

\newtheorem{theorem}{Theorem}
\newtheorem{lemma}{Lemma}[section]
\newtheorem{observation}[lemma]{Observation}
\crefname{observation}{Observation}{Observations}
\Crefname{observation}{Observation}{Observations}
\newtheorem{proposition}{Proposition}[section]

\newtheorem{numberedclaim}{Claim}[theorem]

\newenvironment{claimproof}
  {\begingroup\begin{proof}}
  {\end{proof}\endgroup}

\counterwithin{lemma}{section}
\counterwithout{proposition}{section}

\crefname{algocf}{algorithm}{algorithms}
\Crefname{algocf}{Algorithm}{Algorithms}

\definecolor{purple}{HTML}{B34A20}
\DeclareEmphSequence{\color{purple}\itshape}

\AtBeginEnvironment{algorithm}{%
  \DeclareEmphSequence{\itshape}
}
\AtEndEnvironment{algorithm}{%
  \DeclareEmphSequence{\color{purple}\itshape}
}

\AtBeginEnvironment{thebibliography}{%
  \DeclareEmphSequence{\itshape}
}
\AtEndEnvironment{thebibliography}{%
  \DeclareEmphSequence{\color{purple}\itshape}
}

\DeclareMathOperator{\diam}{diam}

\renewcommand{\int}{\operatorname{int}}

\newcommand{\calB}{\mathcal{B}}
\newcommand{\calE}{\mathcal{E}}
\newcommand{\calI}{\mathcal{I}}
\newcommand{\calM}{\mathcal{M}}

\title{A Low-Diameter Nash-Williams--Tutte Theorem}
\author{
Pin-Hsian Lee\thanks{National Taiwan University} \and 
Peng-Ting Lin\thanks{Academia Sinica} \and
Meng-Tsung Tsai\footnotemark[2] \and
Ting-An Wu\footnotemark[2]
}
\date{}

\begin{document}

\maketitle

\begin{abstract}
The Nash-Williams--Tutte theorem characterizes the existence of $k$ edge-disjoint spanning trees by a condition on every vertex partition. We prove that, for simple graphs, the same condition guarantees $k$ edge-disjoint spanning trees, each of diameter $O(n/\delta)$ and with $O((n/\delta)\log(\delta+1))$ internal vertices, where $n$ is the number of vertices and $\delta$ is the minimum degree. Both bounds hold simultaneously for every tree, with absolute implicit constants independent of $k$, and each is asymptotically optimal in general. Our proof is constructive and yields a randomized algorithm with expected running time $\tilde{O}(n^5)$.

\medskip

We derive applications to short edge-disjoint routes and weak linkages, small Steiner trees, and survivable networks with short paths. The joint bounds also yield low-diameter spanning Eulerian subgraphs and short tours for postman and general routing problems.
\end{abstract}

\thispagestyle{empty}
\clearpage
\pagenumbering{arabic}
\setcounter{page}{1}

\section{Introduction}\label{sec:intro}

The theorem of \citet{Nash-Williams61} and \citet{Tutte61} gives an exact condition for a graph to contain \(k\) edge-disjoint spanning trees: every partition of its vertex set into \(r\) nonempty sets must have at least \(k(r-1)\) edges joining distinct sets. Such trees provide a common representation of edge-disjoint routes between every pair of vertices. The classical theorem, however, does not control the lengths of these routes or the number of internal vertices in each tree. We show that, for simple graphs, the same partition condition suffices to retain all \(k\) trees while simultaneously bounding both parameters. The bounds are asymptotically optimal and have absolute implicit constants, independent of \(k\).

\begin{theorem}[Low-diameter Nash-Williams--Tutte theorem]\label{thm:main}
Let \(G\) be an \(n\)-vertex simple graph with \(n\ge2\) and minimum degree \(\delta\), and let \(k\ge1\) be an integer. The graph \(G\) contains \(k\) edge-disjoint spanning trees if and only if every partition of its vertex set into \(r\) nonempty sets has at least \(k(r-1)\) edges joining distinct sets. Whenever this condition holds, the trees can be chosen so that each has diameter \(O(n/\delta)\) and \(O((n/\delta)\log(\delta+1))\) internal vertices, simultaneously. All implicit constants are absolute. Moreover, these trees can be computed by a randomized algorithm in expected \(\tilde{O}(n^5)\) time.
\end{theorem}

This running-time bound follows from a direct algorithmic implementation of our proof.

Since the partition condition implies \(\delta\ge k\), the theorem also gives diameter \(O(n/k)\) and \(O((n/k)\log(k+1))\) internal vertices in every tree. The two bounds in terms of \(\delta\) are separately asymptotically optimal in general, already for one spanning tree: the extremal examples of \citet*{ErdosPPT89} force diameter \(\Omega(n/\delta)\), while connected-domination lower bounds force \(\Omega((n/\delta)\log(\delta+1))\) internal vertices; see \citet{Alon2023}. Simplicity is essential. Replacing each edge of an \(n\)-vertex path by \(k\) parallel edges gives a multigraph containing \(k\) edge-disjoint spanning trees, but every spanning tree has diameter \(n-1\) and \(n-2\) internal vertices.

\paragraph{Low-diameter spanning trees.}
For a single spanning tree, the \(O(n/\delta)\) diameter bound follows from the classical radius bounds of \citet*{ErdosPPT89}. \citet{Hasunuma00} proved that a simple graph containing \(k\) edge-disjoint spanning trees contains \(k\) such trees, each of diameter \(O(kn/\delta)\). This gives the desired order for fixed \(k\), but the constant grows with \(k\). \citet*{ChandraCDGL24} gave a distributed construction of \(\Omega(\lambda/\log n)\) edge-disjoint spanning trees of diameter \(O((n/\delta)\log n)\), where \(\lambda\) is the edge connectivity. In a different parameter regime, \citet*{ChuzhoyPT20} obtained \(\Omega(\lambda/\log n)\) edge-disjoint spanning trees of diameter \(O(\lambda^{D(D+1)/2})\) when \(\lambda\) is sufficiently large compared with \(\log n\), where \(D\) is the diameter of the host graph. These results use different parameters and, in the distributed setting, a different computational model. Our result gives a centralized construction retaining every tree permitted by the exact partition condition, with a constant in the \(O(n/\delta)\) bound that is independent of \(k\).

\paragraph{Spanning trees with many leaves.}
The number of internal vertices in a spanning tree is closely related to connected domination. Combining the domination bound of \citet{Lovasz1975} with the connected-domination bound of \citet{DuchetMeyniel1982} gives $O((n/\delta)\log(\delta+1))$ internal vertices in a single spanning tree, with asymptotic leading constant $3$ when natural logarithms are used. \citet{KleitmanWest1991} reduced this constant to arbitrarily close to $2.5$, and \citet*{CaroWestYuster2000} reduced it to $1$. \citet{Alon2023} subsequently sharpened the bound to $\frac{n}{\delta+1}(\ln(\delta+1)+O(1))$. For graphs containing $k$ edge-disjoint spanning trees, the degree-balancing results of \citet*{IllingworthPowierskiScottTamitegama2023} and \citet{Horsch2023}, together with a balanced distribution of the remaining edges, yield $k$ edge-disjoint connected spanning subgraphs of minimum degree $\Omega_k(\delta)$. Applying the single-tree bounds within these subgraphs gives $k$ edge-disjoint spanning trees, each with $O_k((n/\delta)\log(\delta+1))$ internal vertices. \Cref{thm:main} removes the dependence on $k$.

\paragraph{Simultaneous diameter and internal-vertex bounds.}
For a single spanning tree, both bounds follow by combining existing results. \citet*{ChandranDasRajendraprasadVarma2012} give a connected two-step dominating set $S$ of size $O(n/\delta)$, while the classical domination bound gives a dominating set $D$ of size $O((n/\delta)\log(\delta+1))$. Attaching $D$ to a spanning tree of the subgraph $G[S]$ induced by $S$ through a forest of depth at most two, and then attaching all remaining vertices as leaves adjacent to $D$, gives a tree of diameter at most $|S|+5$ with $O(|S|+|D|)$ internal vertices. The same decomposition into connected spanning subgraphs used above yields $k$ edge-disjoint spanning trees, each simultaneously having diameter $O_k(n/\delta)$ and $O_k((n/\delta)\log(\delta+1))$ internal vertices. \Cref{thm:main} removes the dependence on $k$.

\begin{graytheorembox}
\Cref{thm:main} retains all \(k\) edge-disjoint spanning trees and simultaneously achieves both optimal bounds on every tree, with absolute implicit constants independent of \(k\).
\end{graytheorembox}

For the routing applications below, assume that \(G\) satisfies the partition condition for \(k\), and let \(\delta\) denote its minimum degree.

\paragraph{Short routes and prescribed terminal pairs.}
\citet{GalilY95} proved a short version of Menger's theorem: in a simple graph of minimum degree at least \(k\), any two vertices joined by \(k\) edge-disjoint paths are joined by \(k\) such paths, each of length \(O(n/k)\). A degree-sensitive refinement of the proof of \citet*{HenzingerKR97} gives length \(O(n/\delta)\) when \(\delta\ge k\). In our setting, \cref{thm:main} gives a fixed tree representation for all vertex pairs simultaneously. The trees use exactly \(k(n-1)\) edges, and their unique paths supply the \(k\) routes for any queried pair. 

The same representation preserves arbitrary prescribed pairings. 
Let $T_1,\ldots,T_k$ be the trees supplied by \cref{thm:main}. For prescribed pairs $(s_i,t_i)$, pick the unique path in $T_i$ from $s_i$ to $t_i$. The resulting paths are edge-disjoint and each has length \(O(n/\delta)\), even when terminals repeat. Thus every \(2k\)-edge-connected simple graph admits a short version of \emph{weak \(k\)-linkage}, which requires edge-disjoint paths joining any prescribed \(k\) terminal pairs. \citet{Huck1991} established weak \(k\)-linkage under \((k+1)\)-edge-connectivity for odd \(k\) and \((k+2)\)-edge-connectivity for even \(k\). Our consequence adds length control under 
a stronger connectivity hypothesis for larger \(k\).

\paragraph{Small Steiner trees and short survivable networks.}
The internal-vertex bound additionally controls the size of the subtrees needed to connect terminal sets. For any set of \(s\ge1\) terminals, pruning each tree to its minimal subtree containing them gives \(k\) edge-disjoint Steiner trees, each of diameter \(O(n/\delta)\) and with \(O(s+(n/\delta)\min\{s,\log(\delta+1)\})\) edges. The minimum combines two bounds: a connector is the union of \(s-1\) short tree paths, and every nonterminal vertex it retains is internal in the original tree.

This observation also applies to \emph{unweighted survivable network design}, where each unordered pair of distinct vertices \(u,v\) has a nonnegative integer requirement \(r_{uv}\) for the number of edge-disjoint connecting paths. Let \(q=\max_{u,v}r_{uv}\le k\), and let \(\mathrm{OPT}\) be the minimum number of edges in a subgraph meeting these requirements without length restrictions. We construct a subgraph with at most \(2\mathrm{OPT}+O((qn/\delta)\log(\delta+1))\) edges in which every pair \(u,v\) has \(r_{uv}\) edge-disjoint paths, each of length \(O(n/\delta)\). At level \(i\), we prune a different tree to the vertices incident with a requirement of at least \(i\). Any feasible subgraph has degree at least \(\max_u r_{uv}\) at vertex \(v\), so the total number of these vertices, counted over all levels, is at most \(2\mathrm{OPT}\). The internal-vertex bound controls the additional vertices. The guarantee complements the classical cost approximation of \citet{Jain2001}: it bounds route lengths while allowing an explicit additive term relative to the unconstrained optimum.

The same connectors, combined with paths that pair up the odd-degree vertices, yield \(k\) alternative feasible postman and general-routing tours with pairwise disjoint augmentation edge sets. The connector bounds control the additional traversals; the required edges may be shared by the tours.

\begin{graytheorembox}
These applications use both guarantees on the same trees: the diameter bound controls route lengths, while the internal-vertex bound controls connector size and branching. Together they yield small Steiner trees and survivable networks with short routes, as well as low-diameter Eulerian subgraphs with few branching vertices.
\end{graytheorembox}

\paragraph{Further structural consequences.}
The theorem also gives quantitative refinements of classical graph structures. Let \(G\) be \(\lambda\)-edge-connected, with \(\lambda\ge4\).
\begin{itemize}
\item \textit{Spanning Eulerian subgraphs with low diameter and few branching vertices.}
\citet{Jaeger1979} showed that every graph containing two edge-disjoint spanning trees is \textit{supereulerian}: it contains a connected spanning Eulerian subgraph. Applying this construction to pairs of trees from \cref{thm:main} yields \(\lfloor\lambda/4\rfloor\) edge-disjoint connected spanning Eulerian subgraphs, each of diameter \(O(n/\delta)\), with at most \(2n-2\) edges and only \(O((n/\delta)\log(\delta+1))\) vertices of degree greater than two. All branching vertices lie among the internal vertices of the two trees used in the construction.

\item \textit{Short cycles through prescribed edges.}
The spanning-tree deletion theorem, as presented by \citet*{CatlinLaiShao2009}, guarantees \(\lfloor\lambda/2\rfloor\) edge-disjoint spanning trees after deleting any \(\lfloor\lambda/2\rfloor\) edges. The remaining minimum degree is at least \(\delta/2\). Hence \cref{thm:main} places any prescribed \(\lfloor\lambda/2\rfloor\) distinct edges on that many edge-disjoint cycles, each of length \(O(n/\delta)\) and containing exactly one prescribed edge.

\item \textit{Shallow spanning out-branchings.}
In the setting where the orientation may be chosen, rather than fixed as in \citet{Edmonds1973}, directing the trees away from any prescribed root gives \(\lfloor\lambda/2\rfloor\) arc-disjoint spanning out-branchings, each of depth \(O(n/\delta)\) and with \(O((n/\delta)\log(\delta+1))\) internal vertices.

\item \textit{Low-diameter spanning trees sharing a common bipartition.}
The maximum-cut argument of \citet*{BangJensenHavetKriesellYeo2022} gives a bipartite spanning subgraph with edge-connectivity at least \(\lceil\lambda/2\rceil\) and minimum degree at least \(\delta/2\). Applying \cref{thm:main} to this subgraph yields \(\lfloor(\lambda+1)/4\rfloor\) edge-disjoint spanning trees sharing a common bipartition, each satisfying both the diameter and internal-vertex bounds.
\end{itemize}

\subsection{Proof Overview of the Main Theorem}\label{sec:overview}

We outline the construction in \cref{thm:main}. Given an \(n\)-vertex
simple graph \(G\) containing \(k\) edge-disjoint spanning trees, we
first find such trees \(T_1,\ldots,T_k\) in \(\tilde{O}(n^3)\) time
using the algorithm of \citet{GabowWestermann1992} or
\citet{BlikstadMukhopadhyayNanongkaiTu2023}. For sufficiently large
\(k\), we transform these trees in Steps 1--4 and then strengthen the
result in Step 5. For bounded \(k\), we let \(T'''_i=T_i\) for every
\(i\in[k] \coloneqq \{1, 2, \ldots, k\}\) and proceed directly to Step 5.

\begin{description}
\item[Step 1: Supply leaves.]
Repartition the edges of \(T_1,\ldots,T_k\) into spanning trees
\(T'_1,\ldots,T'_k\) so that, outside a set \(S\) of size \(O(n/k)\),
every vertex is a leaf in at least \(k/128\) of the first
\(\lceil k/2\rceil\) trees. This follows from
\cref{thm:leaf-supplying} and takes expected
\(\tilde{O}(k^{7/2}n^{3/2})\) time.

\item[Step 2: Group the trees.]
Let \(J=\{1,\ldots,\lceil k/2\rceil\}\) and \(I=[k]\setminus J\).
The trees \(T'_j\) with \(j\in J\) supply leaf edges for shortening
the trees \(T'_i\) with \(i\in I\).

\item[Step 3: Shorten the second group.]
Apply \cref{thm:improvement} to \(T'_1,\ldots,T'_k\), using \(J\) as
the supplying class and \(I\) as the receiving class. The resulting
trees \(T''_1,\ldots,T''_k\) are edge-disjoint, and every \(T''_i\)
with \(i\in I\) has diameter \(O(n/k)\) and
\(O((n/k)\log k)\) internal vertices. Moreover, outside a set \(S'\)
of size \(O(n/k)\), every vertex is a leaf in at least \(k/4\) trees
in this group. This step takes expected \(O(kn^3\log k)\) time.

\item[Step 4: Shorten the first group.]
Apply \cref{thm:improvement} to \(T''_1,\ldots,T''_k\), using \(I\)
as the supplying class and \(J\) as the receiving class. In the
resulting trees \(T'''_1,\ldots,T'''_k\), every tree with an index
in \(J\) satisfies the two bounds. The changes to the trees with
indices in \(I\) preserve their bounds. Thus all \(k\) trees are
edge-disjoint, each with diameter \(O(n/k)\) and
\(O((n/k)\log k)\) internal vertices. This step takes expected
\(O(kn^3\log k)\) time.

\item[Step 5: Replace \(k\) by \(\delta\).]
If \(\delta=1\), the tree \(T'''_1\) already satisfies the bounds
in \cref{thm:main}. Otherwise, apply \cref{cor:minimum-degree} to
\(G\). The resulting trees have diameter
\(O(n/\delta)\) and \(O((n/\delta)\log(\delta+1))\) internal
vertices. This step takes expected \(\tilde{O}(n^5)\) time.
\end{description}
In total, the construction takes expected \(\tilde{O}(n^5)\) time.
The expected bounds include the time for rejected samples and their
retries.

\subsection{Hardness of Improving the Bounds}

The diameter and internal-vertex bounds in \cref{thm:main} are asymptotically optimal. Moreover, even for graphs promised to contain $k$ edge-disjoint spanning trees, deciding whether all $k$ trees have diameter at most $c(n/\delta)$ is NP-hard for a sufficiently small constant $c>0$ and every fixed $k\ge2$ (\cref{thm:diameter-hardness}). Similarly, deciding whether all $k$ trees have at most $c(n/\delta)\log(\delta+1)$ internal vertices is NP-hard for a sufficiently small constant $c>0$ and every sufficiently large fixed $k$ (\cref{thm:internal-vertices-hardness}).

\paragraph{Paper Organization.}
\cref{sec:prel} introduces notation and preliminaries. \cref{sec:leaf-supplying,sec:improvement,sec:outdegree-shortening} establish the tools used in \cref{sec:combine} to prove \cref{thm:main}. \cref{sec:hardness-results} proves the hardness results, and \cref{app:set-cover-hardness} proves the Set Cover hardness theorem used in \cref{sec:hardness-results}.
\section{Notation and Preliminaries} 
\label{sec:prel}

For a positive integer $k$, let $[k]=\{1,\ldots,k\}$. For a function $f$ on a finite set and a subset $F$ of its domain, let $f(F)=\sum_{e\in F}f(e)$.

For a graph $G$, let $V(G)$ and $E(G)$ denote its vertex set and edge set, respectively. We write $uv$ for an undirected edge $\{u,v\}$. A collection $\mathcal{H}$ of subgraphs of $G$ is \emph{edge-disjoint} if $E(H)\cap E(H')=\emptyset$ for any two distinct $H,H'\in\mathcal{H}$.

For a vertex $v\in V(G)$, let $N_G(v)$ denote its neighborhood and let $\deg_G(v)=|N_G(v)|$ denote its degree. We write $\delta(G)$ and $\lambda(G)$ for the minimum degree and edge connectivity of $G$, respectively. For $S\subseteq V(G)$, let $G[S]$ denote the subgraph of $G$ induced by $S$. For vertices $s,t$ of a connected graph $G$, let $\operatorname{dist}_G(s,t)$ denote the length of a shortest $s$--$t$ path, where length is measured in edges. The diameter $\operatorname{diam}(G)$ is the maximum distance between two vertices of $G$.

For a tree $T$, let $\int(T)=\{v\in V(T):\deg_T(v)\ge2\}$ denote its set of internal vertices. A vertex of degree one in a tree is a \emph{leaf}.

An \emph{orientation} of a graph assigns a direction to each edge. For an arc $(u,v)$, the vertex $u$ is its \emph{tail}, $v$ is its \emph{head}, and $v$ is an \emph{outneighbor} of $u$. The \emph{outdegree} of a vertex is the number of arcs directed away from it.

A family $\mathcal{F}$ of sets is \emph{laminar} if, for every $A,B\in\mathcal{F}$, either $A\subseteq B$, $B\subseteq A$, or $A\cap B=\emptyset$. We follow standard notation for matroid theory; see, for example, \citet{Oxley2011}.

Unless otherwise stated, $\log$ denotes the natural logarithm. The notation $\tilde{O}$ suppresses polylogarithmic factors in the relevant parameters, and a subscript on asymptotic notation indicates parameters on which the implicit constants may depend.

\section{Supplying Leaves for Shortening}
\label{sec:leaf-supplying}

\begin{theorem}\label{thm:leaf-supplying}
Let $G=(V,E)$ be an $n$-vertex simple graph with $n\ge2$ whose edge set is the union of $k\ge1$ edge-disjoint spanning trees. There is an edge partition of $G$ into spanning trees $T_1,\ldots,T_k$ such that every vertex outside a set $S$ of size $O(n/k)$ is a leaf in at least $k/128$ of $T_1,\ldots,T_{\lceil k/2\rceil}$.
\end{theorem}

The theorem supplies many leaf occurrences at almost every vertex while retaining all $k$ trees. We obtain the partition by repeatedly extracting a tree that has few internal vertices in a randomly selected set. The key requirement is that each extraction leave a spanning-tree partition of the remaining edges.

We state the extraction lemma (\cref{lem:leafy-extraction}) and use it to prove the theorem. The extraction lemma is proved afterward using a structural bound on tight sets (\cref{lem:tight-cut-core}), a weighted estimate (\cref{lem:weighted-bound}), and the weight-splitting theorem of \citet{Frank1981}.

\subsection{
From one tree to many leaf occurrences}

\begin{lemma}\label{lem:leafy-extraction}
Let $G=(V,E)$ be an $n$-vertex simple graph with $n\ge2$ whose edge set is the union of $k\ge1$ edge-disjoint spanning trees. Let $L\subseteq V$ and $r\in V$, and orient each of these trees toward $r$. Suppose that every vertex has at most $k/8$ outneighbors in $L$. Let $T$ minimize $\sum_{e\in E(T)}w(e)$ among all spanning trees of $G$ for which $E\setminus E(T)$ can be partitioned into $k-1$ spanning trees, where $w(e)=|e\cap L|$ for every $e\in E$. Then $|\int(T)\cap L|=O(n/k)$.
\end{lemma}

The orientation condition makes $L$ sparse among outgoing neighbors. The lemma ensures that minimizing the number of incidences with $L$ produces a tree in which all but $O(n/k)$ vertices of $L$ are leaves, without preventing further extractions.

\begin{proof}[Proof of Theorem~\ref{thm:leaf-supplying}]
If $k<16$, let $S=V$ and choose any edge partition of $G$ into $k$ spanning trees. Assume $k\ge16$, and let $m=\lceil k/2\rceil$. Fix a root $r$. We extract $T_1,\ldots,T_m$ successively, maintaining a spanning-tree partition of the unused edges.

At round $i$, let $t=k-i+1\ge k/2$ be the number of remaining trees, and orient them toward $r$. Sample each vertex independently with probability $1/16$ to form $Y_i$, using randomness independent of earlier rounds. Call a nonroot vertex bad if at least $t/8$ of its outneighbors lie in $Y_i$. Let $L_i$ be obtained from $Y_i$ by deleting every sampled outneighbor of every bad vertex. Every vertex then has at most $t/8$ outneighbors in $L_i$. Apply Lemma~\ref{lem:leafy-extraction} to extract a tree $T_i$ with $|\int(T_i)\cap L_i|=O(n/t)$, and retain a partition of the remaining edges into $t-1$ trees.

Each nonroot vertex has $t$ distinct outneighbors, so a Chernoff bound gives a bad-vertex probability of at most $e^{-t/48}$. Each bad vertex causes at most $t$ deletions. Consequently,
\[
\begin{aligned}
\mathbb{E}[|Y_i\cap\int(T_i)|]
&\le \mathbb{E}[|Y_i\setminus L_i|]+O(n/t)\\
&\le nt e^{-t/48}+O(n/t)
=O(n/k).
\end{aligned}
\]
Here $t e^{-t/48}=O(1/t)$ and $t\ge k/2$. Let $F$ count all sampled internal-vertex occurrences over the $m$ rounds. Then
\[
F=\sum_{i=1}^{m}|Y_i\cap\int(T_i)|,
\qquad
\mathbb{E}[F]=O(n).
\]

Let $Z_v=|\{i:v\in Y_i\}|$, and let $S$ consist of the vertices that are leaves in fewer than $k/128$ extracted trees. The independent samples give $Z_v\sim\operatorname{Bin}(m,1/16)$, even though the extracted trees depend on earlier rounds. Since $\mathbb{E}[Z_v]=m/16\ge k/32$, another Chernoff bound gives
\[
\Pr[Z_v<k/64]\le e^{-k/256}.
\]
Every vertex $v$ in $S$ with $Z_v\ge k/64$ contributes at least $k/128$ to $F$. Therefore,
\[
|S|\le |\{v:Z_v<k/64\}|+\frac{128F}{k},
\qquad
\mathbb{E}[|S|]\le ne^{-k/256}+\frac{128}{k}\mathbb{E}[F]=O(n/k).
\]
Some choice of the samples therefore gives $|S|=O(n/k)$. Retain its extracted trees and label the remaining trees $T_{m+1},\ldots,T_k$. Every vertex outside $S$ is a leaf in at least $k/128$ of the first $m$ trees.
\end{proof}

\subsection{Controlling tight sets}

The following lemma depends only on the given tree partition and its orientation, not on $L$ or any edge weights. Its role is to isolate vertices that can account for most of the arcs leaving a tight set.

\begin{lemma}\label{lem:tight-cut-core}
Let $G=(V,E)$ be an $n$-vertex simple graph with $n\ge2$ whose edges are partitioned into $k\ge1$ spanning trees, each oriented toward $r\in V$. Call a nonempty set $A\subseteq V$ \emph{tight} if $|E(G[A])|=k(|A|-1)$, and call it \emph{nontrivial} if $|A|\ge2$. Let $R$ consist of $r$ and every vertex $v$ for which some nontrivial tight set $A$ containing $v$ has at least $3k/4$ arcs directed from $v$ to $V\setminus A$. Then $|R|=O(n/k)$.
\end{lemma}
\begin{proof}
For any tight set $A$, each given tree contributes at most $|A|-1$ edges to $G[A]$, so equality forces every tree to induce a tree on $A$. Every vertex other than $r$ has $k$ outgoing arcs. Subtracting the $k(|A|-1)$ arcs internal to $A$ shows that exactly $k$ arcs leave $A$ if $r\notin A$, and none leave $A$ otherwise.

We record the effect of intersecting two tight sets $A,B$. Let $h$ count the edges joining $A\setminus B$ to $B\setminus A$. Since every nonempty set $S$ satisfies $|E(G[S])|\le k(|S|-1)$,
\[
\begin{aligned}
k(|A|+|B|-2)+h
&=|E(G[A\cup B])|+|E(G[A\cap B])|\\
&\le k(|A\cup B|-1)+k(|A\cap B|-1)\\
&=k(|A|+|B|-2).
\end{aligned}
\]
Thus $h=0$, and both $A\cup B$ and $A\cap B$ are tight.

For each $v\in R\setminus\{r\}$, let $A_v$ be an inclusion-minimal witnessing set. Passing to a smaller nontrivial tight set containing $v$ preserves the witnessing arcs, so $A_v$ is inclusion-minimal among all nontrivial tight sets containing $v$. In particular, $r\notin A_v$, and distinct vertices of $R\setminus\{r\}$ have distinct sets $A_v$, since two such vertices would supply more than $k$ outgoing arcs from the same tight set.

Suppose that $|A_v\cap A_w|\ge2$ and neither set contains the other. The intersection is a nontrivial tight set, so minimality implies $v\in A_v\setminus A_w$ and $w\in A_w\setminus A_v$. Since no edge joins these exclusive parts, the witnessing arcs from both vertices leave the tight union $A_v\cup A_w$, contradicting its outgoing-arc bound. Hence the sets $A_v$ are nested or intersect in at most one vertex.

If $A_v\subseteq A_w$ and $A_v\ne A_w$, minimality gives $w\notin A_v$. At least $3k/4$ of the $k$ arcs leaving $A_w$ start at $w$, so at most $k/4$ start at $v$. At least $3k/4$ arcs from $v$ leave $A_v$, so at least $k/2$ end in $A_w\setminus A_v$. By simplicity,
\[
|A_w\setminus A_v|\ge k/2.
\]

Let $T$ be one of the given spanning trees. Since each $T[A_v]$ is a nontrivial tree, its edge set determines $A_v$. Thus the sets $E(T[A_v])$ are distinct and form a laminar family on the $n-1$ edges of $T$: sets whose vertex sets intersect in at most one vertex have disjoint edge sets. Moreover,
\[
k(|A_v|-1)=|E(G[A_v])|\le\binom{|A_v|}{2},
\]
so $|A_v|\ge2k$ and each member has at least $2k-1$ edges. For nested members,
\[
|E(T[A_w])\setminus E(T[A_v])|
=|A_w\setminus A_v|\ge k/2.
\]
Form the inclusion forest of this family, with each nonmaximal member a child of its smallest strict superset in the family. For each leaf, take its associated edge set. For each node with exactly one child, take the difference between its associated edge set and that of its child. These sets are pairwise disjoint, and each contains $\Omega(k)$ edges. Since they are subsets of $E(T)$, there are $O(n/k)$ leaves and nodes with exactly one child. The number of branching nodes is at most the number of leaves, so $|R\setminus\{r\}|=O(n/k)$. Since $n\ge2k$ by simplicity, adding $r$ preserves this bound.
\end{proof}

\subsection{Bounding spanning-tree weights}

Using Lemma~\ref{lem:tight-cut-core}, we establish the weighted inequality in \cref{lem:weighted-bound}, which we then use to prove Lemma~\ref{lem:leafy-extraction}.

\begin{lemma}\label{lem:weighted-bound}
Let $G=(V,E)$ be an $n$-vertex simple graph with $n\ge2$ whose edge set is the union of $k\ge8$ edge-disjoint spanning trees. Let $L\subseteq V$ and $r\in V$, and orient each of these trees toward $r$. Suppose that every vertex has at most $k/8$ outneighbors in $L$. Define $w(e)=|e\cap L|$ for every $e\in E$, and let $\mu(f)$ denote the maximum $f$-weight of a spanning tree of $G$. For every $p:E\to\mathbb Z$, let $p(E)=\sum_{e\in E}p(e)$. Then
\[
p(E)-(k-1)\mu(p)-\mu(p-w)\le |L|+O(n/k).
\]
\end{lemma}
\begin{proof}
Let $R$ be the set supplied by Lemma~\ref{lem:tight-cut-core}, let $X=V\setminus L$, and let $g(p)=p(E)-(k-1)\mu(p)-\mu(p-w)$. Adding an integer constant to every $p(e)$ leaves $g(p)$ unchanged because $|E|=k(n-1)$. We may therefore assume that $p(e)-w(e)\ge1$ for every $e \in E$.

For each integer $j\ge0$, let
\[
G_j=(V,\{e\in E:p(e)\ge j\}),
\qquad
H_j=(V,\{e\in E:p(e)-w(e)\ge j\}).
\]
Then $G_{j+2}\subseteq H_j\subseteq G_j$. Let $\mathcal{S}_j$ be the collection of vertex sets of the connected components of $G_j$. For $C\in\mathcal{S}_j$, let $q_j(C)$ count the components of $H_j[C]$, and let $\sigma_j(C)=k(|C|-1)-|E(G_j[C])|\ge0$. Since $H_j\subseteq G_j$, each component of $H_j[C]$ is also a component of $H_j$. Kruskal's algorithm expresses $\mu(p)$ and $\mu(p-w)$ as sums of threshold ranks. Within $C$, these ranks are $|C|-1$ and $|C|-q_j(C)$, respectively. Counting $p(E)$ by thresholds $j$ therefore gives
\[
\begin{aligned}
g(p)
&=\sum_{j\ge1}\sum_{C\in\mathcal{S}_j}
\bigl(|E(G_j[C])|-(k-1)(|C|-1)-(|C|-q_j(C))\bigr)\\
&=\sum_{j\ge1}\sum_{C\in\mathcal{S}_j}
\bigl(q_j(C)-1-\sigma_j(C)\bigr).
\end{aligned}
\]
Singleton sets $C$ contribute zero.

Call a vertex set $U$ \emph{large} if $|U|\ge k/8$ or $U\cap R\ne\emptyset$. For $C\in\mathcal{S}_j$, let $b_j(C)$ count the large components of $H_j[C]$, and let $t_j(C)$ count the components $U\subseteq L$ joined to $X$ by an edge of $G_j[C]$. These two collections need not be disjoint. Our goal is to prove the component bound
\[
q_j(C)-1-\sigma_j(C)
\le t_j(C)+\max\{b_j(C)-1,0\},
\]
and then sum its two terms by charging vertices of $L$ and counting branches of a splitting tree.

For $|C|\ge2$, we show that at most $\sigma_j(C)$ components are counted by neither $b_j(C)$ nor $t_j(C)$. For each such component $U$, let $v\in U\cap X$ if possible, and otherwise let $v\in U$. No edge of $G_j[C]$ joins $v$ to $X\setminus U$: when $v\in X$, its $w$-weight would be zero and it would also lie in $H_j[C]$; otherwise $U\subseteq L$ and the edge would make $U$ count toward $t_j(C)$. Since $U\cap R=\emptyset$ and $r\in R$, $v$ has $k$ outgoing arcs. At most $k/8$ lead to $L$ by hypothesis, and fewer than $k/8$ lead to $U$ because $|U|<k/8$. Thus $v$ supplies at least $3k/4$ outgoing arcs whose underlying edges are absent from $G_j[C]$.

There are at most $k|C|-|E(G_j[C])|=k+\sigma_j(C)$ arcs whose tails lie in $C$ and whose underlying edges are absent from $G_j[C]$. If more than $\sigma_j(C)$ components were of neither type, their distinct representatives would therefore give
\[
(3k/4)(\sigma_j(C)+1)
\le k|C|-|E(G_j[C])|
=k+\sigma_j(C).
\]
Hence $\sigma_j(C)\le k/(3k-4)<1$, so $\sigma_j(C)=0$ by integrality. The spanning-tree partition then gives
\[
|E(G_j[C])|=k(|C|-1)\ge |E(G[C])|,
\]
so $G_j[C]=G[C]$ and $C$ is tight. Each chosen representative therefore has at least $3k/4$ outgoing arcs leaving $C$, contradicting its being outside $R$. Consequently, $q_j(C)\le b_j(C)+t_j(C)+\sigma_j(C)$, proving the component bound. For $|C|=1$, we have $q_j(C)-1-\sigma_j(C)=0$, so the component bound also holds.

It remains to bound the two terms on the right-hand side.

For every component $U$ counted by $t_j(C)$, let $\{v,x\}\in E(G_j[C])$ with $v\in U$ and $x\in X$, and charge $U$ to $v$. Because $v\in U\subseteq L$, $x\in X$, and $w(\{v,x\})\coloneqq|\{v,x\}\cap L|$, we have $w(\{v,x\})=1$. This edge is absent from $H_j$, since $U$ is a component of $H_j$ disjoint from $X$. Together with its membership in $G_j$, this gives $p(\{v,x\})=j$ by integrality. An edge $\{v,y\}\in E$ with $y\in X$ and larger $p$-weight would have $w(\{v,y\})=1$ and belong to $H_j$, contradicting that $U$ is a component disjoint from $X$. Hence
\[
j=\max\{p(\{v,y\}):y\in X,\ \{v,y\}\in E\}.
\]
Thus $v$ determines its only possible charging threshold and belongs to only one component of $H_j$ there. Each vertex of $L$ receives at most one charge, giving
\[
\sum_j\sum_{C\in\mathcal{S}_j}t_j(C)\le |L|.
\]

For the branching term, consider the two sequences
\[
G_0,H_0,G_2,H_2,\ldots
\qquad\text{and}\qquad
G_1,H_1,G_3,H_3,\ldots.
\]
Along either sequence, components only split, and the first graph is $G$ by normalization. For each sequence, include both $G_j$ and $H_j$ at its first index $j>\max_{e\in E}p(e)$ and stop there. Let $S$ be the \emph{splitting tree} with one node for every large component at every level of the sequence. The node associated with $V$ at the first level is the root. Each nonroot node in $S$ is a child of the unique node at the preceding level whose associated component contains its vertex set. This node exists because every ancestor of a large component is large.

The leaves of $S$ have disjoint vertex sets. At most $8n/k$ have size at least $k/8$, and at most $|R|$ others contain a vertex of $R$. A node in $S$ corresponding to $C\in\mathcal{S}_j$ has $b_j(C)$ children; a component not represented in $S$ has $b_j(C)=0$. The sum of the number of children minus one over all branching nodes in $S$ is the number of leaves minus one. The branching term uses only the $G_j$ levels; including the $H_j$ levels can only increase it. Summing over both trees therefore gives
\[
\sum_j\sum_{C\in\mathcal{S}_j}\max\{b_j(C)-1,0\}
=O(n/k+|R|)=O(n/k).
\]
Summing the component bound and applying these two estimates proves the lemma.
\end{proof}
\subsection{Proof of the extraction lemma}

For an edge set $F$, let $w(F)=\sum_{e\in F}w(e)$. Since
\[
w(E(T))=\sum_{v\in L}\deg_T(v)\ge |L|+|\int(T)\cap L|,
\]
it suffices to bound the minimum feasible tree weight by $|L|+O(n/k)$. Frank's weight-splitting theorem will express this weight in the form bounded by Lemma~\ref{lem:weighted-bound}.

\begin{theorem}[Weight-splitting theorem;~\citet{Frank1981}]
\label{thm:frank}
Let $M_1$ and $M_2$ be matroids on the same ground set $E$, with respective families of bases $\mathcal{B}_1$ and $\mathcal{B}_2$, and suppose that $\mathcal{B}_1\cap\mathcal{B}_2\ne\emptyset$. For every integer-valued weight function $w:E\to\mathbb{Z}$, there exist a basis $B \in \mathcal{B}_1 \cap \mathcal{B}_2$ and integer-valued weight functions $w_1,w_2:E\to\mathbb{Z}$ with $w=w_1+w_2$ such that $B$ is a maximum $w_1$-weight basis in $\mathcal{B}_1$ and a maximum $w_2$-weight basis in $\mathcal{B}_2$.
\end{theorem}

\begin{proof}[Proof of Lemma~\ref{lem:leafy-extraction}]
If $k<8$, then $|\int(T)\cap L|\le n\le7n/k$, so assume $k\ge8$. Let $\calE=\{(e,i):e\in E,\ i\in[k]\}$. Let $\calE_e=\{(e,i):i\in[k]\}$ for every $e\in E$ and $\calE_i=\{(e,i):e\in E\}$ for every $i\in[k]$. Within each $\calE_i$, identify $(e,i)$ with $e$ when referring to a subgraph of $G$ or its weight. Define $\calM_1=(\calE,\calI_1)$ and $\calM_2=(\calE,\calI_2)$ by
\begin{itemize}
\item $\calI_1=\{F\subseteq\calE:\forall i\in[k],\ F\cap\calE_i\mbox{ is a forest in }G\}$, and
\item $\calI_2=\{F\subseteq\calE:\forall e\in E,\ |F\cap\calE_e|\le1\}$.
\end{itemize}

$\calM_1$ is the direct sum of $k$ graphic matroids, and $\calM_2$ is a partition matroid. Let $\calB_1$ and $\calB_2$ be their respective collections of bases. An edge-partition of $G$ into $k$ spanning trees gives a common independent set of cardinality $k(n-1)$, the rank of both matroids. Hence $\calB_1\cap\calB_2\ne\emptyset$. Conversely, every common basis induces an edge-partition of $G$ into $k$ spanning trees.

Define $c:\calE\to\mathbb Z$ by $c(e,1)=-w(e)$ and $c(e,i)=0$ for $i\ge2$. Every common basis $B$ satisfies $c(B)=-w(B\cap\calE_1)$. Thus maximizing $c(B)$ minimizes $w(B\cap\calE_1)$ subject to the remaining edges being partitionable into $k-1$ spanning trees. By the choice of $T$, the maximum common-basis weight is $-w(E(T))$.

By Theorem~\ref{thm:frank}, there are integer-valued functions $c_1,c_2:\calE\to\mathbb Z$ with $c=c_1+c_2$ and a common basis $B$ that is a maximum $c_1$-weight basis of $\calM_1$ and a maximum $c_2$-weight basis of $\calM_2$. Such a basis also maximizes $c$ among all common bases.

We replace $c_2(e,i)$ by $c'_2(e,i):=q(e)$ and $c_1(e,i)$ by $c'_1(e,i):=c(e,i)-q(e)$, where $q(e)=\max_{i\in[k]}c_2(e,i)$. These replacements preserve $c'_1+c'_2=c$ and the optimality of $B$ in both matroids. Indeed, every basis of $\calM_2$ has the same $c'_2$-weight. Moreover, the original optimality of $B$ in $\calM_2$ implies that $c_2(e,i)=q(e)$ for every $(e,i)\in B$. Thus $c'_1\le c_1$, with equality on $B$, so $B$ remains optimal in $\calM_1$.

Let $p(e)=-q(e)$, and let $\mu(f)$ denote the maximum $f$-weight of a spanning tree of $G$. By the optimality of $B$, we have $c'_1(B)=\mu(p-w)+(k-1)\mu(p)$ and $c'_2(B)=-p(E)$. Since $c'_1(B)+c'_2(B)=c(B)=-w(E(T))$, Lemma~\ref{lem:weighted-bound} gives
\[
w(E(T))=p(E)-(k-1)\mu(p)-\mu(p-w)\le |L|+O(n/k).
\]

Finally, $w(E(T))=\sum_{v\in L}\deg_T(v)\ge |L|+|\int(T)\cap L|$, so $|\int(T)\cap L|=O(n/k)$.
\end{proof}

\subsection{Running time}
\label{subsec:leaf-supplying-algorithm}

The construction in \cref{thm:leaf-supplying} can be implemented as a Las Vegas algorithm with expected running time \(\tilde{O}(k^{7/2}n^{3/2})\), and hence \(\tilde{O}(n^5)\). The input consists of \(G\), represented by adjacency lists, and \(k\); an initial spanning-tree partition is not required. The algorithm returns the entire tree partition and the exceptional set \(S\).  Under the theorem's hypotheses, \(|E|=k(n-1)\) and \(k\le n/2\).

\paragraph{Implementing an extraction.}
Suppose that \(t\) trees remain. The two matroids in the proof of Lemma~\ref{lem:leafy-extraction}, with \(k\) replaced by \(t\), have ground-set size \(N_t=t^2(n-1)\) and common rank \(\rho_t=t(n-1)\). Both have explicit graphic representations. For \(\calM_1\), use \(t\) disjoint copies of the residual graph. For \(\calM_2\), use a separate pair of vertices for each residual edge \(e\), joined by \(t\) parallel edges corresponding to \((e,1),\ldots,(e,t)\). A forest in this second graph contains at most one copy of each \(e\). Both graphs can be constructed in \(O(N_t+tn)\) time.

The extraction requires a maximum-weight common basis, rather than an unrestricted maximum-weight common independent set with the nonpositive weights \(c\). To enforce the required cardinality, give each element weight \(a(e,i)=2n-1+c(e,i)\). A common basis has weight at least \((2n-1)\rho_t-2(n-1)\), since only its first tree incurs a cost. Any smaller common independent set has weight at most \((2n-1)(\rho_t-1)\). Thus a maximum-weight common independent set is a common basis and, among common bases, maximizes \(c\). Its first tree is the required minimum-weight extraction, and its other trees explicitly partition the residual edges.

The weighted graphic matroid-intersection algorithm of \citet{GabowXu1989} runs in \[
O(N\sqrt{\rho}\log^2(2\rho)\log(2\rho W))
\]
time for two explicitly represented graphic matroids on \(N\) elements of rank at most \(\rho\), with positive integer weights at most \(W\); see also \citet*{HuangKakimuraKamiyama2019}. Here \(W\le2n-1\), so one extraction takes \(O(t^{5/2}n^{3/2}\log^3(kn))\) time. Frank's weight splitting, the set \(R\), and the threshold graphs are used only to prove the extraction guarantee; the algorithm need not construct them.

\paragraph{Obtaining and processing the tree partitions.}
We compute an initial partition once, using the same two graphic representations with \(t=k\) and uniform weight one. Maximizing the weight then maximizes cardinality, and the input promise guarantees a common basis of size \(k(n-1)\). This costs \(O(k^{5/2}n^{3/2}\log^3(kn))\). If \(k<16\), return this partition and \(S=V\), as in the existence proof.

For \(k\ge16\), one attempt performs the \(\lceil k/2\rceil\) extraction rounds from that proof. At a round with \(t\) remaining trees, rooting and orienting these trees, sampling \(Y_i\), counting sampled outneighbors, and marking the vertices deleted to obtain \(L_i\) take \(O(tn)\) time. Constructing the graphic representations and their weights is already included in the extraction bound. Summing over the rounds gives \(O(k^{7/2}n^{3/2}\log^3(kn))\) time per attempt. Computing all leaf counts and the resulting set \(S\), verifying the returned tree partition, and writing the output take \(O(kn)\) additional time. In particular, the bound includes reading the \(k(n-1)\) input edges.

\paragraph{Verifying success and bounding retries.}
Let \(C\) be a fixed absolute constant for which the proof of \cref{thm:leaf-supplying} gives \(\mathbb{E}[|S|]\le Cn/k\). After an attempt, accept precisely when \(|S|\le2Cn/k\); this test uses the actual leaf counts. Markov's inequality gives success probability at least \(1/2\). Upon rejection, restart from the saved initial partition with fresh random samples. All other choices can be made deterministically. The expectation bound remains valid despite the dependence of each extracted tree on earlier rounds, because the proof conditions on those rounds before applying the outneighbor estimate. Hence the expected number of attempts is at most two. Every accepted output satisfies the required leaf guarantee, and every attempt has the time bound above, including unsuccessful attempts. Together with the initialization cost, this proves the claimed expected running time.
\section{Shortening Theorem}
\label{sec:improvement}

In this section, we prove a shortening theorem that repartitions the edges between two classes of spanning trees. The trees in one class attain the desired bounds on diameter and the number of internal vertices, while those in the other class incur only additive increases in these quantities.

\begin{theorem}[Shortening theorem]\label{thm:improvement}\label{lem:main}
Let \(k\) be sufficiently large. Let \(G=(V,E)\) be a simple graph on \(n\ge2\) vertices whose edges are partitioned into spanning trees \(T_1,\ldots,T_k\). Partition \([k]\) into a receiving class \(I\) and a supplying class \(J\). Suppose that every vertex outside a set \(S\subseteq V\) is a leaf in at least \(k/128\) supplying trees.

Then \(E\) has a new partition into spanning trees \(T_1',\ldots,T_k'\). For each receiving tree \(T_i'\), \(i\in I\), we have
\[
\diam(T_i')=O\!\left(|S|+\frac nk\right),
\qquad
|\int(T_i')|=O\!\left(|S|+\frac nk\log k\right).
\]
For each supplying tree \(T_j'\), \(j\in J\), we have
\[
\diam(T_j')\le\diam(T_j)+2,
\qquad
|\int(T_j')|\le|\int(T_j)|+O\!\left(\frac nk\log k\right).
\]
There is also a set \(S'\supseteq S\) of size \(O(|S|+n/k)\) outside which each vertex is internal in fewer than \(k/8\) receiving trees. Given the original partition, the construction takes \(O(kn^3\log k)\) expected time.
\end{theorem}

We call the trees indexed by \(I\) the \emph{receiving trees} and those indexed by \(J\) the \emph{supplying trees}. In particular, when \(|S|=O(n/k)\), both the diameters of the receiving trees and \(|S'|\) are \(O(n/k)\). The proof starts by reserving leaf edges from the supplying trees (\cref{lem:restoration-matchings}). Each such edge belongs to the reserve at its leaf endpoint, which we call its \emph{owner}. Whenever a receiving tree takes a reserved edge at \(v\), it returns another edge incident with \(v\) (\cref{lem:shorten-receiving}). Thus \(v\) still has the same number of edges to use when we reattach it to the supplying trees, even though their other endpoints may have changed. The proof is constructive, and we analyze the expected running time of the resulting algorithm in \cref{sec:lemma2_runtime}.

\subsection{Shortening one receiving tree}

The following lemma uses reserved edges from the supplying trees to shorten a receiving tree and bound its number of internal vertices, while preserving the size of each reserve.

\begin{lemma}[Shortening one receiving tree]\label{lem:shorten-receiving}
Let \(T\) be a spanning tree of a simple graph \(G=(V,E)\) on \(n\ge2\) vertices. For each \(v\in V\), let \(R_v\subseteq E\) be a set of edges incident with \(v\). Assume that \(E(T)\) and the sets \((R_v)_{v\in V}\) are pairwise disjoint. Let \(A,B\subseteq V\) and let \(d\ge2\) be an integer. Suppose that every \(v\in V\setminus B\) has at least \(d\) edges in \(R_v\) whose other endpoint belongs to \(A\).

The set \(E(T)\cup\bigcup_{v\in V}R_v\) can be partitioned into the edge set of a spanning tree \(T'\) and sets \((R_v')_{v\in V}\), where each \(R_v'\) consists of edges incident with \(v\) and \(|R_v'|=|R_v|\). The tree satisfies
\[
\diam(T')=O\!\left(|B|+\frac{|A|}{d}\right),
\qquad
|\int(T')\setminus(A\cup B)|\le2\left(|B|+\frac{|A|}{d}\right).
\]
\end{lemma}

The proof uses the following lemma to shorten a tree while preserving the size of each reserve. We prove this auxiliary lemma in \cref{sec:outdegree-shortening}.

\begin{lemma}[Outdegree-preserving shortening]\label{lem:shorten}
Let \(T\) be a spanning tree of a simple graph \(G=(V,E)\) on \(N\ge1\) vertices. Let \(X\subseteq V\) and let $d$ be a positive integer. Suppose that \(D\) is an orientation of \(E\setminus E(T)\) in which every vertex outside \(X\) has outdegree at least \(d\). Then there are a spanning tree \(Q\) of \(G\) and an orientation \(D'\) of \(E\setminus E(Q)\) such that
\[
\diam(Q)=O(|X|+N/d),
\]
and every vertex has the same outdegree in \(D'\) as in \(D\). Such a tree and orientation can be found deterministically in \(O(dN^3)\) time.
\end{lemma}

\begin{proof}[Proof of \cref{lem:shorten-receiving}]
\emph{Build a rooted forest.}
Start with one rooted component for each vertex of \(A\cup B\). Whenever a root \(v\notin B\) has an edge of \(R_v\) to a vertex of \(A\) in a different component, add that edge as the parent edge of \(v\). The root of the other component becomes the root of the merged tree. Let \(F_0\) be the forest when no further merge is possible.

Every vertex of \(B\) remains a root, so each component that meets \(B\) contains exactly one vertex of \(B\). In every other component, the root belongs to \(A\setminus B\), and all its reserved neighbors in \(A\) lie in the same component. There are at least \(d\) distinct neighbors because \(G\) is simple. Counting the root as well, each such component contains at least \(d+1\) vertices of \(A\setminus B\). Thus the number \(c\) of components satisfies
\[
c\le|B|+\frac{|A\setminus B|}{d+1}.
\]
Now attach each vertex outside \(A\cup B\) to a vertex of \(A\) by one edge of its own reserve. The resulting spanning forest \(F\) still has \(c\) components. Every vertex outside \(A\cup B\) is a leaf of \(F\), and every nonroot vertex \(v\) owns its parent edge.

\emph{Insert the forest into the tree.}
Starting from \(T\), insert the edges of \(F\) in order from the roots toward the leaves. To insert the parent edge of \(v\), take it from \(R_v\) and add it to the current tree. On the resulting cycle, remove the other edge incident with \(v\) and return it to \(R_v\). This exchange preserves the tree and replaces one edge of \(R_v\) by another incident with \(v\).

The edge removed at \(v\) lies outside \(F\). Indeed, \(F\) was disjoint from \(T\), and none of the child edges at \(v\) has yet been inserted. Thus the final tree \(U\) contains \(F\), and every reserve retains its size. \Cref{fig:insertion} illustrates two exchanges.

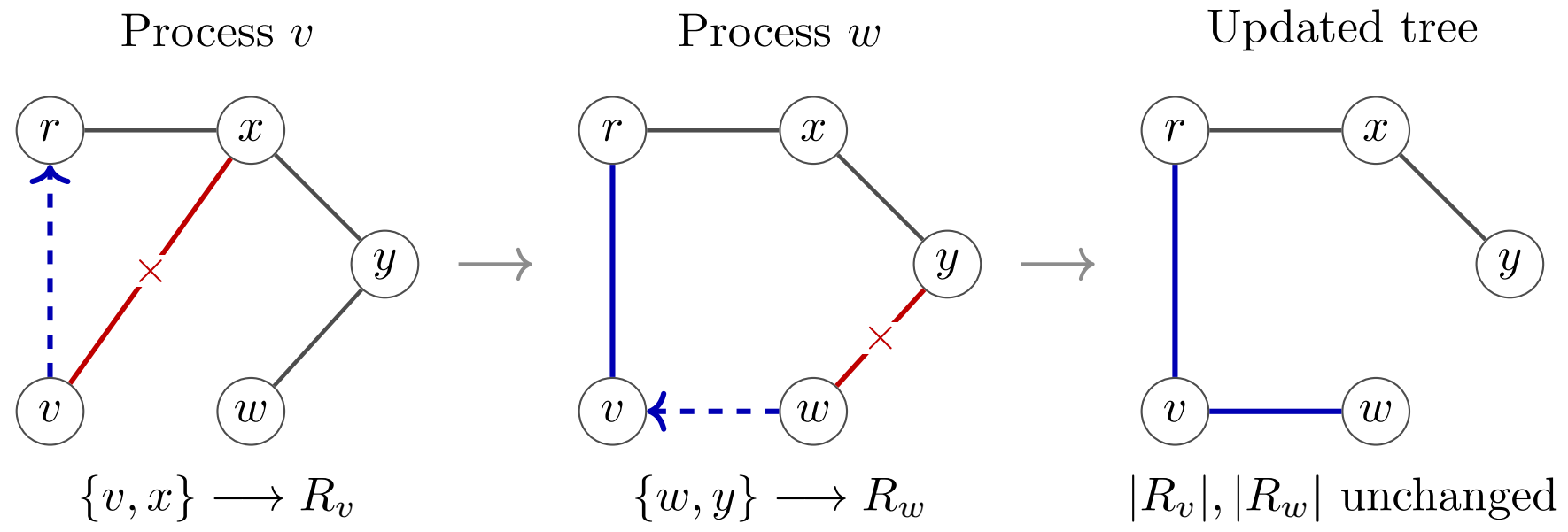
\begin{figure}[H]
\centering
\begin{tikzpicture}[x=1cm,y=1cm,font=\small,
 vertex/.style={circle,draw=black!70,fill=white,inner sep=2pt,minimum size=5mm},
 treeedge/.style={black!70,line width=0.85pt},
 inserted/.style={blue!70!black,line width=1.15pt},
 inserting/.style={blue!70!black,line width=1.15pt,dashed,->},
 removed/.style={red!75!black,line width=1.05pt}]
\begin{scope}
\node at (1.25,2.65) {Process \(v\)};
\node[vertex] (r) at (0,1.9) {\(r\)};
\node[vertex] (x) at (1.5,1.9) {\(x\)};
\node[vertex] (y) at (2.5,0.9) {\(y\)};
\node[vertex] (v) at (0,-0.2) {\(v\)};
\node[vertex] (w) at (1.5,-0.2) {\(w\)};
\draw[treeedge] (r)--(x)--(y)--(w);
\draw[removed] (v)--(x);
\node[text=red!75!black,fill=white,inner sep=0pt] at (0.75,0.85) {\(\times\)};
\draw[inserting] (v)--(r);
\node at (1.25,-0.85) {\(\{v,x\}\longrightarrow R_v\)};
\end{scope}
\draw[->,black!45,line width=0.8pt] (3.05,0.9)--(3.6,0.9);
\begin{scope}[xshift=4.2cm]
\node at (1.25,2.65) {Process \(w\)};
\node[vertex] (r) at (0,1.9) {\(r\)};
\node[vertex] (x) at (1.5,1.9) {\(x\)};
\node[vertex] (y) at (2.5,0.9) {\(y\)};
\node[vertex] (v) at (0,-0.2) {\(v\)};
\node[vertex] (w) at (1.5,-0.2) {\(w\)};
\draw[treeedge] (r)--(x)--(y);
\draw[removed] (w)--(y);
\node[text=red!75!black,fill=white,inner sep=0pt] at (2,0.35) {\(\times\)};
\draw[inserted] (v)--(r);
\draw[inserting] (w)--(v);
\node at (1.25,-0.85) {\(\{w,y\}\longrightarrow R_w\)};
\end{scope}
\draw[->,black!45,line width=0.8pt] (7.25,0.9)--(7.8,0.9);
\begin{scope}[xshift=8.4cm]
\node at (1.25,2.65) {Updated tree};
\node[vertex] (r) at (0,1.9) {\(r\)};
\node[vertex] (x) at (1.5,1.9) {\(x\)};
\node[vertex] (y) at (2.5,0.9) {\(y\)};
\node[vertex] (v) at (0,-0.2) {\(v\)};
\node[vertex] (w) at (1.5,-0.2) {\(w\)};
\draw[treeedge] (r)--(x)--(y);
\draw[inserted] (r)--(v)--(w);
\node at (1.25,-0.85) {\(|R_v|,|R_w|\) unchanged};
\end{scope}
\end{tikzpicture}
\caption{Inserting the path \(w\to v\to r\). First \(v\) takes \(vr\) from its reserve and returns \(vx\) to it. Next \(w\) adds \(wv\) and returns \(wy\) to its own reserve. The second exchange leaves \(vr\) in the tree.}
\label{fig:insertion}
\end{figure}

\emph{Identify the core.}
The tree \(U\) has exactly \(c-1\) edges outside \(F\), joining the \(c\) components of \(F\). A leaf of \(F\) can become internal in \(U\) only at an endpoint of one of these joining edges. Let \(D\) be the set of their endpoints. Then
\[
|D|\le2(c-1)\le2\left(|B|+\frac{|A|}{d}\right).
\]
Set \(W=A\cup B\cup D\). Every vertex outside \(W\) is still a leaf attached to \(A\). Deleting these leaves gives the \emph{core}, the tree \(U[W]\) induced by \(W\) (\cref{fig:core}). Its order is \(N:=|W|\le|A|+3|B|+2|A|/d\).

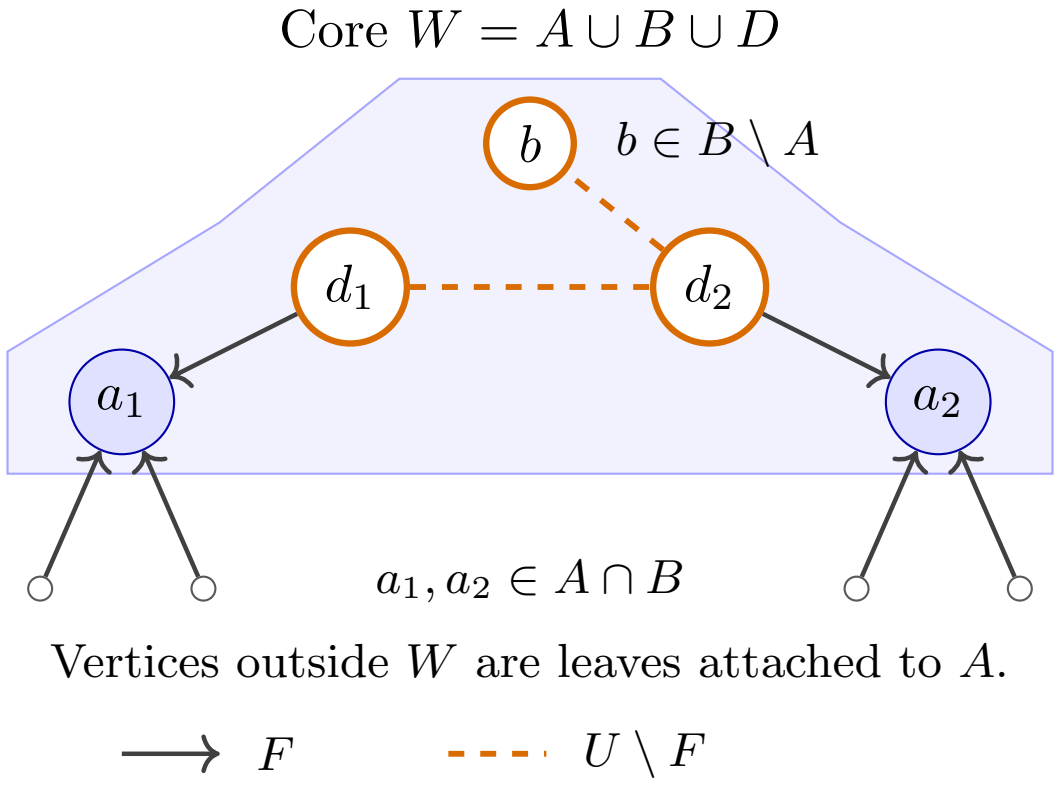
\begin{figure}[!htbp]
\centering
\begin{tikzpicture}[x=1.08cm,y=0.95cm,font=\small,
 fedge/.style={draw=black!75,line width=0.85pt,->},
 joinedge/.style={draw=orange!85!black,line width=1.1pt,dashed},
 anode/.style={circle,draw=blue!65!black,fill=blue!12,inner sep=3pt},
 dnode/.style={circle,draw=orange!85!black,line width=1.2pt,fill=white,inner sep=3pt},
 leaf/.style={circle,draw=black!65,fill=white,minimum size=4.5pt,inner sep=0pt}]
\filldraw[fill=blue!5,draw=blue!35,thin] (-0.7,0.15)--(-0.7,1)--(0.6,1.9)--(1.7,2.9)--(3.3,2.9)--(4.4,1.9)--(5.7,1)--(5.7,0.15)--cycle;
\node at (2.5,3.25) {Core \(W=A\cup B\cup D\)};
\node[anode] (a1) at (0,0.65) {\(a_1\)};
\node[anode] (a2) at (5,0.65) {\(a_2\)};
\node[dnode] (x) at (1.4,1.45) {\(d_1\)};
\node[dnode] (y) at (3.6,1.45) {\(d_2\)};
\node[dnode] (b) at (2.5,2.45) {\(b\)};
\node[anchor=west,font=\footnotesize] at (2.9,2.45) {\(b\in B\setminus A\)};
\node[leaf] (l1) at (-0.5,-0.65) {};
\node[leaf] (l2) at (0.5,-0.65) {};
\node[leaf] (l3) at (4.5,-0.65) {};
\node[leaf] (l4) at (5.5,-0.65) {};
\draw[fedge] (x)--(a1);
\draw[fedge] (y)--(a2);
\draw[fedge] (l1)--(a1);
\draw[fedge] (l2)--(a1);
\draw[fedge] (l3)--(a2);
\draw[fedge] (l4)--(a2);
\draw[joinedge] (x)--(y);
\draw[joinedge] (y)--(b);
\node[font=\footnotesize] at (2.5,-0.6) {\(a_1,a_2\in A\cap B\)};
\node[font=\footnotesize] at (2.5,-1.15) {Vertices outside \(W\) are leaves attached to \(A\).};
\draw[fedge] (0,-1.8)--(0.6,-1.8);
\node[anchor=west,font=\footnotesize] at (0.7,-1.8) {\(F\)};
\draw[joinedge] (2,-1.8)--(2.6,-1.8);
\node[anchor=west,font=\footnotesize] at (2.7,-1.8) {\(U\setminus F\)};
\end{tikzpicture}
\caption{The dashed edges join the three components of \(F\). Their endpoints form \(D=\{b,d_1,d_2\}\). Only the shaded core \(W\) is shortened. The solid arrows point from children to parents in \(F\). In this example \(a_1\) and \(a_2\) belong to both \(A\) and \(B\).}
\label{fig:core}
\end{figure}

\emph{Shorten the core.}
Each vertex used at most one of its original reserved edges, namely its parent edge in \(F\). Hence each \(v\in W\setminus B\) still has at least \(d-1\) reserved edges connecting to \(A\). Select \(d-1\) of them and orient them away from \(v\). Both ends of every selected edge lie in \(W\). Apply \cref{lem:shorten} to \(U[W]\) together with the selected edges, taking \(B\) as the exceptional set. We obtain a tree \(Q\) on \(W\) and an orientation of the edges left over in this graph. Since \(d-1\ge d/2\), the bound on \(N\) gives
\[
\diam(Q)=O\!\left(|B|+\frac{N}{d-1}\right)
=O\!\left(|B|+\frac{|A|}{d}\right).
\]
Replace \(U[W]\) by \(Q\) and keep the leaf edges outside \(W\) in place. Remove the selected edges from the reserves and assign each complementary edge to its tail. By \cref{lem:shorten}, each vertex receives exactly as many edges as it gave up, yielding the required sets \(R_v'\).

A path in the resulting tree \(T'\) uses at most one pendant edge at each end, so \(\diam(T')\le\diam(Q)+2\). If \(B\) is empty, \(A\) contains a root of \(F_0\) and its \(d\) distinct reserved neighbors. Hence \(|B|+|A|/d\ge1\), which absorbs the additive two. Every internal vertex of \(T'\) lies in \(W=A\cup B\cup D\), giving the second bound.
\end{proof}

\subsection{Proof of the shortening theorem}

We now prove \cref{thm:improvement}. The next lemma lets us choose the reserves before the first exchange and restore the supplying trees after the last. We prove it in \cref{sec:restoration-matchings}.

\begin{lemma}[Restoring the supplying trees]\label{lem:restoration-matchings}
Let \(k\) be sufficiently large. Let \(G=(V,E)\) be a simple graph on \(n\ge2\) vertices whose edges are partitioned into spanning trees \(T_1,\ldots,T_k\). Let \(I,J\) partition \([k]\). Suppose every vertex outside \(S\subseteq V\) is a leaf in at least \(k/128\) of the trees indexed by \(J\). Set \(r=\lfloor k/8192\rfloor\).

There exist a set \(H\subseteq V\) and pairwise disjoint edge sets \((R_v)_{v\in V}\) satisfying the following.
\begin{enumerate}
\item[(i)] \(H\supseteq S\) and \(|H\setminus S|\le n/k\).

\item[(ii)] Each \(R_v\) consists of edges incident with \(v\) from the supplying trees. Its size is \(r\) for \(v\notin H\), and \(R_v=\emptyset\) for \(v\in H\).

\item[(iii)] Let
\[
P=\bigcup_{i\in I}E(T_i)\cup\bigcup_{v\in V}R_v
\]
be the set of edges available for exchange. Suppose \(P\) is repartitioned into spanning trees \((T_i')_{i\in I}\) and reserves \((R_v')_{v\in V}\), each \(R_v'\) consisting of \(|R_v|\) edges incident with \(v\). Then the edges of \(G\) outside the new receiving trees can be partitioned into spanning trees \((T_j')_{j\in J}\) such that, for every \(j\in J\),
\[
\diam(T_j')\le\diam(T_j)+2,
\qquad
|\int(T_j')|\le|\int(T_j)|+O\!\left(\frac nk\log k\right).
\]
\end{enumerate}
\end{lemma}

\begin{proof}[Proof of \cref{thm:improvement}]
Take \(H\) and the reserves from \cref{lem:restoration-matchings}. Write \(R_v\) for the current reserve at \(v\). Its size stays at \(r=\lfloor k/8192\rfloor\) for \(v\notin H\).

\emph{Choose attachment points for each tree.}
Choose \(d=\lceil\log k\rceil\) and \(b=\lceil96(n/r)d\rceil\). For sufficiently large \(k\), we have \(b\le n\). Process the receiving trees in a fixed order. When we reach \(T_i\), choose a uniform \(b\)-element set \(A_i\subseteq V\), independently of all earlier choices. Put a vertex \(v\notin H\) in \(X_i\) if fewer than \(d\) edges of \(R_v\) join it to \(A_i\). Then every vertex outside \(B_i=H\cup X_i\) meets the hypothesis of \cref{lem:shorten-receiving}.

Shorten \(T_i\) using that lemma. Let \(D_i=\int(T_i')\setminus(A_i\cup B_i)\) be the internal vertices of the new tree outside \(A_i\cup B_i\). We have
\begin{equation}\label{eq:tree}
\begin{aligned}
\diam(T_i')&=O(|B_i|+b/d),\qquad |D_i|\le2(|B_i|+b/d),\\
\int(T_i')&\subseteq A_i\cup B_i\cup D_i.
\end{aligned}
\end{equation}

\emph{Count appearances across trees.}
Recall that we need a set \(S'\) of size \(O(|S|+n/k)\) such that every vertex outside \(S'\) is internal in fewer than \(k/8\) receiving trees. Outside \(H\cup\bigcup_iX_i\), every internal vertex of \(T_i'\) belongs to \(A_i\) or \(D_i\). To form \(S'\), add the vertices that appear at least \(k/16\) times in either family. Set
\[
\begin{aligned}
Z_A&=\{v\in V:|\{i\in I:v\in A_i\}|\ge k/16\},\\
Z_D&=\{v\in V:|\{i\in I:v\in D_i\}|\ge k/16\},
\end{aligned}
\qquad
S'=H\cup\bigcup_{i\in I}X_i\cup Z_A\cup Z_D.
\]
For \(v\notin S'\), the two counts add to less than \(k/8\). This proves the required bound on the number of trees in which \(v\) is internal. To bound \(|S'|\), observe that each vertex of \(Z_D\) contributes at least \(k/16\) to \(\sum_{i\in I}|D_i|\). Hence
\begin{equation}\label{eq:exceptional-count}
|S'|\le|H|+\sum_{i\in I}|X_i|+|Z_A|+\frac{16}{k}\sum_{i\in I}|D_i|.
\end{equation}

\emph{Bound the size of the exceptional set.}
First bound \(|X_i|\). Condition on all choices before we sample \(A_i\). This fixes the current reserve at each \(v\notin H\), whose \(r\) edges lead to distinct neighbors because \(G\) is simple. The number \(Y\) of these neighbors in \(A_i\) is hypergeometric with mean \(\mu=rb/n\ge96d\). The Chernoff bound gives
\[
\Pr[v\in X_i\mid\text{earlier choices}]\le e^{-rb/(8n)}\le k^{-10},
\qquad
\mathbb{E}\sum_{i\in I}|X_i|\le nk^{-9}.
\]
For \(Z_A\), we count how often a fixed vertex \(v\) is itself sampled. The sets \(A_i\) are independent, so \(z(v)=|\{i\in I:v\in A_i\}|\) is binomial with mean \(|I|b/n=O(\log k)\). The Chernoff bound gives, for sufficiently large \(k\),
\[
\Pr[z(v)\ge k/16]\le\left(\frac{16e|I|b}{kn}\right)^{k/16}\le k^{-10},
\qquad
\mathbb{E}|Z_A|\le nk^{-10}.
\]
Combining these estimates with \eqref{eq:exceptional-count}, \(|I|\le k\), and \(|D_i|\le2(|H|+|X_i|+b/d)\) gives
\[
\begin{aligned}
\mathbb{E}|S'|&\le33|H|+32b/d+(1+32/k)nk^{-9}+nk^{-10}\\
&\le33(|H|+b/d+n/k).
\end{aligned}
\]
By Markov's inequality, with probability at least \(1/2\) we have
\begin{equation}\label{eq:outcome}
|S'|\le66(|H|+b/d+n/k).
\end{equation}

\emph{Finish the partition.}
Fix a sample satisfying \eqref{eq:outcome}. Since \(|H|\le|S|+n/k\) and \(b/d=O(n/k)\), we have \(|S'|=O(|S|+n/k)\). Every \(B_i\) lies in \(S'\), so \eqref{eq:tree} gives \(\diam(T_i')=O(|S|+n/k)\). Counting the internal vertices in \(A_i\), \(B_i\), and \(D_i\) gives
\[
|\int(T_i')|\le b+|B_i|+2(|B_i|+b/d)
=O\!\left(|S|+\frac nk\log k\right).
\]
Each reserve still has its original size and consists of edges incident with its owner. Thus \cref{lem:restoration-matchings} restores the supplying trees with the stated bounds and completes the partition of \(E\).
\end{proof}

\subsection{Restoring the supplying trees}
\label{sec:restoration-matchings}

We now prove \cref{lem:restoration-matchings}. We select leaf edges from the supplying trees, reserving half for exchanges with the receiving trees and protecting the other half for reconstruction. After the exchanges, the reserved edges may have different endpoints. We use matchings to reattach the deleted leaves to the supplying trees, and restrict the attachment points to small sets to control the number of new internal vertices.

The following lemma gives a degree condition for a perfect matching and shows that it survives sampling with high probability. We prove it at the end of this subsection.
\begin{lemma}[Sampling a perfect matching]\label{lem:sampled-matching}
Let \(B\) be a simple bipartite graph with vertex classes \(L\) and \(R\), both of size \(m\ge2\). Let \(\delta_L\) and \(\delta_R\) be the minimum degrees on \(L\) and \(R\). Suppose that \(\delta_L+\delta_R\ge5m/4\).
\begin{enumerate}
\item[(i)] \(B\) has a perfect matching.
\item[(ii)] Independently for each \(j\in R\), choose a random subset \(S_j\subseteq N_B(j)\). Suppose there is a \(p\in(0,1]\) such that
\[
\Pr[S_j\cap X=\emptyset]\le e^{-p|X|}
\qquad\text{for every }j\in R\text{ and }X\subseteq N_B(j).
\]
Let \(B'\) be the spanning subgraph with \(N_{B'}(j)=S_j\) for each \(j\in R\). If \(pm\ge96\log m\), then \(B'\) has a perfect matching with probability at least \(1-4m^{-11}\).
\end{enumerate}
\end{lemma}

\begin{proof}[Proof of \cref{lem:restoration-matchings}]
\emph{Choose the reserves.}

Since \(G\) is simple, \(k(n-1)=|E(G)|\le n(n-1)/2\), so \(n\ge2k\ge3\). For each vertex \(v\), let \(L_v\subseteq J\) index the supplying trees in which \(v\) is originally a leaf. For \(v\notin S\), we have \(|L_v|\ge k/128\ge2r\).

Independently at each \(v\notin S\), choose a uniform \(2r\)-element subset \(J_v\subseteq L_v\), and put \(J_v=\emptyset\) on \(S\). Fix an ordering of \(J\) and protect the leaf edges corresponding to the first \(r\) indices of \(J_v\). The \emph{head} of an edge incident with \(v\) is its other endpoint; let \(P_v\) be the set of protected heads. Put \(P_v=\emptyset\) on \(S\). 

To restore \(v\) using an edge \(vx\) in \(T_j\), we need \(j\notin J_x\); otherwise \(x\) is also selected for deletion from that tree. Since exchanges may change the reserved endpoints, we must consider every possible head \(x\in N_G(v)\).

Set \(q=\lfloor r/4\rfloor\). At each \(v\notin S\), test
\[
\begin{aligned}
|J_v\cap J_x|&\le q
&&\forall x \in N_G(v),\\
|\{x\in P_v:j\in J_x\}|&\le q
&&\forall j \in J_v.
\end{aligned}
\]
The first test requires that each possible head \(x\) is selected for deletion in at most \(q\) trees indexed by \(J_v\). The second requires that at most \(q\) protected heads at \(v\) are selected for deletion from each tree indexed by \(J_v\). We will use these conditions to guarantee a perfect matching that reattaches \(v\) to the supplying trees after the exchanges. Let \(Z\) be the set of vertices where at least one test fails.

Condition on \(J_v\), which also fixes \(P_v\). For \(x\notin S\cup\{v\}\) and any fixed set \(Q\) of \(h\ge1\) indices, sampling without replacement gives
\[
\Pr[Q\subseteq J_x\mid J_v]\le\left(\frac{2r}{|L_x|}\right)^h\le32^{-h}.
\]
For \(x\in S\), this probability is zero. Conditional on \(J_v\), the choices at distinct protected heads are independent, so \(h\) specified heads all select a given \(j\) with probability at most \(32^{-h}\). A failure of either test has a witness of size \(q+1\) among at most \(2r\) indices or heads. Its conditional probability is therefore at most
\[
\binom{2r}{q+1}32^{-(q+1)}\le(e/4)^{q+1}.
\]
There are at most \(2|E(G)|\) tests of the first type and \(2rn\) of the second type. Since \(q=\Theta(k)\), for sufficiently large \(k\),
\[
\mathbb{E}|Z|\le\bigl(2k(n-1)+2rn\bigr)(e/4)^{q+1}
\le3kn(e/4)^{q+1}\le nk^{-10}.
\]
Choose an outcome with \(|Z|\le n/k\), and let \(H=S\cup Z\). For each \(v\in H\), cancel its selections and set \(J_v=P_v=R_v=\emptyset\). Its selected edges remain in their original supplying trees. For each \(v\notin H\), let \(R_v\) consist of its \(r\) unprotected selected leaf edges. The two inequalities still hold at every such \(v\). We have now fixed \(H\) and all the reserves before any repartition, with \(|H\setminus S|\le n/k\).

No selected edge is selected at both ends, since a tree on at least three vertices has no adjacent leaves. Thus all selected edges have a unique owner, their leaf endpoint, and the protected edges and movable reserves are pairwise disjoint.

\emph{Match the available heads after a repartition.}

For each \(j\in J\), delete from \(T_j\) the vertices \(v\) with \(j\in J_v\), and call the retained subtree \(T_j^-\). We delete only original leaves, so \(T_j^-\) is a nonempty tree containing every original internal vertex. Its vertex set is precisely \(\{x\in V:j\notin J_x\}\).

Fix an arbitrary repartition allowed in the lemma, and let \(R_v'\) be the updated reserve at \(v\). For every \(v\notin H\), it still consists of \(r\) edges incident with \(v\), although their other endpoints may have changed. Let \(Y_v\) consist of the other endpoints of these edges and the protected edges at \(v\). Since \(G\) is simple,
\[
P_v\subseteq Y_v\subseteq N_G(v),\qquad |P_v|=r,\qquad |Y_v|=2r.
\]
Define the bipartite graph \(B_v\) with vertex classes \(Y_v\) and \(J_v\) by
\[
xj\in E(B_v)
\quad\Longleftrightarrow\quad
j\notin J_x.
\]
Thus \(xj\) is admissible exactly when the available edge \(vx\) can attach \(v\) directly to the retained tree \(T_j^-\). 

A perfect matching in \(B_v\) would assign one available edge at \(v\) to each supplying tree from which \(v\) was deleted. Choosing such a matching for every \(v\notin H\) would restore spanning trees by attaching the deleted vertices as leaves to the retained subtrees.

We now verify that these matchings exist. The first inequality gives every head in \(Y_v\) at least \(2r-q\) available tree indices. The second gives every tree index in \(J_v\) at least \(r-q\) available protected heads. Consequently,
\[
\delta_{Y_v}+\delta_{J_v}\ge3r-2q\ge\frac54(2r).
\]
Part~(i) of \cref{lem:sampled-matching} therefore gives a perfect matching in every \(B_v\), for every allowed repartition.

Although this restores the supplying trees, the attachment points may become new internal vertices. To bound this increase, we will restrict the heads used in each tree \(j\) to a small set \(C_j\subseteq V(T_j^-)\).

\emph{Restrict the attachment points.}

Independently for each \(j\in J\), choose a uniform subset \(C_j\subseteq V(T_j^-)\) of size \(\min\{b_s,|V(T_j^-)|\}\), where
\[
b_s=\left\lceil48\frac nr\log(2r)\right\rceil.
\]
For sufficiently large \(k\), we have \(b_s\le n\), since \(r=\Theta(k)\). In each \(B_v\), retain exactly the edges \(xj\) with \(x\in C_j\); see \cref{fig:constrained-restoration}.

For a fixed \(v\notin H\), the retained neighbor sets \(C_j\cap N_{B_v}(j)\) are independent over \(j\in J_v\). Put \(n_j=|V(T_j^-)|\). If \(n_j\ge b_s\), then for every \(X\subseteq N_{B_v}(j)\),
\[
\Pr[C_j\cap X=\emptyset]
=\frac{\binom{n_j-|X|}{b_s}}{\binom{n_j}{b_s}}
\le e^{-b_s|X|/n_j}
\le e^{-b_s|X|/n},
\]
where the numerator is zero if \(n_j-|X|<b_s\). If \(n_j<b_s\), then \(C_j=V(T_j^-)\), and the same upper bound holds. Apply part~(ii) of \cref{lem:sampled-matching} with \(m=2r\) and \(p=b_s/n\). Since \(pm=2rb_s/n\ge96\log(2r)\), the restricted \(B_v\) lacks a perfect matching with probability at most \(4(2r)^{-11}\). Call a vertex \emph{failed} if this happens. Linearity of expectation gives at most \(4n(2r)^{-11}\le n/(2k)\) failed vertices in expectation for sufficiently large \(k\); independence between different vertices is not needed. By Markov's inequality, with probability at least \(1/2\), at most \(n/k\) vertices fail. Fix such a sample and choose a surviving matching at every successful vertex. This sample may depend on the fixed repartition; the original reserves do not.

\begin{figure}[!ht]
\centering
\begin{tikzpicture}[x=1cm,y=0.82cm,font=\small,
  vertex/.style={circle,draw=black!65,fill=white,inner sep=1pt,minimum size=5.7mm},
  treeindex/.style={vertex,draw=blue!60!black,fill=blue!4},
  allowed/.style={black,line width=0.65pt},
  matched/.style={blue!65!black,line width=1.4pt}]
\begin{scope}
\node at (1.2,4.75) {Before restriction};
\node at (0,4.2) {Heads};
\node at (2.4,4.2) {Tree indices};
\filldraw[fill=orange!7,draw=orange!65!black,rounded corners=3pt]
  (-0.34,2.10) rectangle (0.34,3.90);
\node[anchor=east] at (-0.43,3.0) {\(P_v\)};
\foreach \name/\lab/\h in {p1/p_1/3.6,p2/p_2/2.4,x1/x_1/1.2,x2/x_2/0}
  \node[vertex] (\name) at (0,\h) {\(\lab\)};
\foreach \i/\h in {1/3.6,2/2.4,3/1.2,4/0}
  \node[treeindex] (j\i) at (2.4,\h) {\(j_\i\)};
\foreach \u/\w in {p1/j2,p2/j1,p2/j3,p2/j4,x1/j1,x1/j4,x2/j2,x2/j3}
  \draw[allowed] (\u.east)--(\w.west);
\foreach \u/\w in {p1/j1,p2/j2,x1/j3,x2/j4}
  \draw[allowed] (\u.east)--(\w.west);
\end{scope}
\draw[->,>=stealth,black!60,line width=0.9pt] (3.05,1.8)--(4.35,1.8);
\node[align=center,font=\footnotesize] at (3.7,0.9) {sample \(C_j\)\\and restrict};
\begin{scope}[xshift=5.15cm]
\node at (1.2,4.75) {After restriction};
\node at (0,4.2) {Heads};
\node at (2.4,4.2) {Tree indices};
\node at (4.9,4.2) {Sampled heads};
\filldraw[fill=orange!7,draw=orange!65!black,rounded corners=3pt]
  (-0.34,2.10) rectangle (0.34,3.90);
\node[anchor=east] at (-0.43,3.0) {\(P_v\)};
\foreach \name/\lab/\h in {p1/p_1/3.6,p2/p_2/2.4,x1/x_1/1.2,x2/x_2/0}
  \node[vertex] (\name) at (0,\h) {\(\lab\)};
\foreach \i/\h in {1/3.6,2/2.4,3/1.2,4/0}
  \node[treeindex] (j\i) at (2.4,\h) {\(j_\i\)};
\foreach \u/\w in {x1/j1,p1/j2,p2/j3,p2/j4}
  \draw[allowed] (\u.east)--(\w.west);
\foreach \u/\w in {p1/j1,p2/j2,x1/j3,x2/j4}
  \draw[allowed] (\u.east)--(\w.west);
\node[anchor=west,font=\footnotesize] at (2.95,3.6)
  {\(C_{j_1}\cap Y_v=\{p_1,x_1\}\)};
\node[anchor=west,font=\footnotesize] at (2.95,2.4)
  {\(C_{j_2}\cap Y_v=\{p_1,p_2\}\)};
\node[anchor=west,font=\footnotesize] at (2.95,1.2)
  {\(C_{j_3}\cap Y_v=\{p_2,x_1\}\)};
\node[anchor=west,font=\footnotesize] at (2.95,0)
  {\(C_{j_4}\cap Y_v=\{p_2,x_2\}\)};
\end{scope}
\end{tikzpicture}
\caption{A schematic illustration of the matching restriction. The sets on the right list the sampled heads for each tree. An original edge \(xj\) survives exactly when \(x\in C_j\). Shading marks protected heads.}
\label{fig:constrained-restoration}
\end{figure}
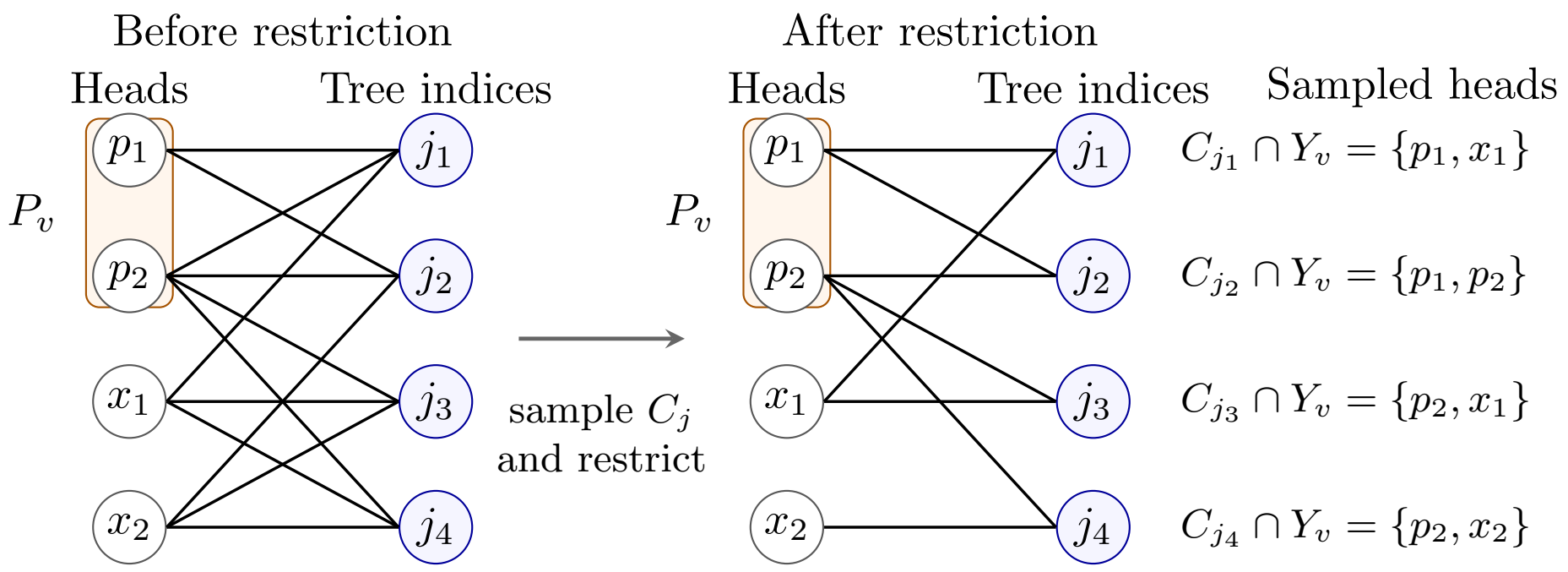

\emph{Repair failures and restore the trees.}
For every failed vertex \(v\), choose a perfect matching in the unrestricted graph \(B_v\) using part~(i) of \cref{lem:sampled-matching}. For each \(j\in J_v\), add its matched head to \(C_j\). This head belongs to \(T_j^-\) by the definition of \(B_v\), and enlarging \(C_j\) preserves all previously chosen matchings. Each failed vertex adds at most one head per tree, so
\[
|C_j|\le b_s+n/k=O\!\left(\frac nk\log k\right)\qquad(j\in J).
\]
The edge sets \(\{vx:x\in Y_v\}\), over \(v\notin H\), contain exactly the protected edges and updated reserves, each with a unique owner. Together with the edges of the retained supplying subtrees and updated receiving trees, they partition \(E(G)\). The matchings assign every available edge to exactly one supplying tree. In the resulting \(T_j'\), every deleted vertex has exactly one edge to \(T_j^-\), so \(T_j'\) is a spanning tree. A longest path in \(T_j'\) consists of a path in \(T_j^-\), possibly a single vertex, and at most one new leaf edge at each end. Only heads in \(C_j\) can become new internal vertices. Hence
\[
\diam(T_j')\le\diam(T_j^-)+2\le\diam(T_j)+2,
\qquad
\int(T_j')\subseteq\int(T_j)\cup C_j.
\]
The bound on \(|C_j|\) gives the required bound on internal vertices. Since the repartition was arbitrary, the reserves fixed above have the stated restoration property for every allowed repartition.
\end{proof}

\begin{proof}[Proof of \cref{lem:sampled-matching}]
By Hall's theorem, a perfect matching fails to exist exactly when some \(U\subseteq L\) has fewer than \(|U|\) neighbors. Equivalently, there is a set \(W\subseteq R\) of size \(m+1-|U|\) with no edge to \(U\).

Fix \(U\subseteq L\) and \(W\subseteq R\) with \(u=|U|\), \(w=|W|\) and \(u+w=m+1\). Let \(e\) count the edges of \(B\) between them. A vertex of \(U\) has at most \(u-1\) neighbors outside \(W\), and a vertex of \(W\) has at most \(w-1\) neighbors outside \(U\). Consequently,
\[
\frac eu+\frac ew\ge(\delta_L-u+1)+(\delta_R-w+1)\ge\frac m4,
\qquad
e\ge\frac{muw}{4(u+w)}\ge\frac m8\min\{u,w\}.
\]
Thus every such pair has an edge, proving (i). For (ii), the pair loses all its edges only if each \(j\in W\) misses its neighbors in \(U\). Independence over \(j\) gives
\[
\begin{aligned}
\Pr[\text{no retained edge joins }U\text{ to }W]
&=\prod_{j\in W}\Pr[S_j\cap(N_B(j)\cap U)=\emptyset]\\
&\le\prod_{j\in W}e^{-p|N_B(j)\cap U|}=e^{-pe}.
\end{aligned}
\]
For a fixed \(s=\min\{u,w\}\), specifying a pair requires choosing the side of the smaller set, its \(s\) vertices, and the \(s-1\) vertices outside the larger set on the opposite side. There are at most \(2m^{2s-1}\) choices. A union bound over these pairs, using \(pm\ge96\log m\), gives
\[
\Pr[B'\text{ has no perfect matching}]
\le\sum_{s\ge1}2m^{2s-1}e^{-pms/8}
\le\frac2m\sum_{s\ge1}m^{-10s}\le4m^{-11}.
\]
This proves (ii).
\end{proof}

\subsection{Running time}
\label{sec:lemma2_runtime}

The construction is a Las Vegas algorithm with expected running time \(O(kn^3\log k)\). We assume that the original trees are given by their edge lists and count operations on \(O(\log n)\)-bit indices. Simplicity gives \(k\le n/2\).

\paragraph{Preparing the reserves.}
Sampling the selections and testing the two conditions take \(O(nk^2)\) time per trial. Each intersection test takes \(O(k)\) time, there are \(O(nk)\) ordered adjacent pairs, and testing the protected heads takes \(O(k^2)\) time per vertex. Repeat until \(|Z|\le n/k\). Since \(\mathbb{E}|Z|\le nk^{-10}\), the expected number of trials is \(O(1)\).

\paragraph{Shortening the receiving trees.}
Each vertex of \(W\setminus B_i\) has complementary outdegree \(d-1\). The implementation of \cref{lem:shorten} analyzed in \cref{sec:outdegree-shortening} therefore takes \(O(dN^3)=O(n^3\log k)\) time on the core, since \(N=|W|\le n\) and \(d=\lceil\log k\rceil\).

The forest construction makes at most \(n-1\) merges, each costing \(O(n^2)\), while the remaining exchanges and updates take \(O(n^2)\). Thus shortening one receiving tree takes \(O(n^3\log k)\) time, and a complete receiving trial takes \(O(kn^3\log k)\). This bound holds for every sampling outcome.

Repeat the receiving trial with fresh samples until \eqref{eq:outcome} holds, resetting the receiving trees and reserves to the accepted preparation after each rejection. Conditional on that preparation, each trial succeeds with probability at least \(1/2\), so the expected number of trials is \(O(1)\).

\paragraph{Restoring the supplying trees.}
Sampling all \(C_j\) takes \(O(nk)\) time. Each matching graph has \(O(k)\) vertices and \(O(k^2)\) edges, so an augmenting-path algorithm tests it in \(O(k^3)\) time. Testing all vertices and repairing the failures therefore take \(O(nk^3)\) time. If more than \(n/k\) vertices fail, resample only the \(C_j\). Conditional on the fixed repartition, each trial succeeds with probability at least \(1/2\), so the expected number of trials is \(O(1)\).

The total expected running time is
\[
O(nk^2+kn^3\log k+nk^3)
=O(kn^3\log k),
\]
using \(k\le n/2\).

\section{Outdegree-Preserving Shortening}
\label{sec:outdegree-shortening}
In this section, we prove \cref{lem:shorten}.

We shorten a spanning tree of large diameter through a sequence of edge exchanges with its complementary orientation, preserving every vertex's outdegree in that orientation. Each exchange decreases a potential defined by the cuts obtained by deleting tree edges. Along a diametral path, the edges whose cuts have sufficient weight on both sides form a middle segment; the remaining edges have a small side. We first bound the total size of these small sides, then the length of the middle segment, and finally the number of edges near the ends. The corresponding algorithm is given after the proof.

\begin{proof}
If \(X=V\), the statement is immediate. Henceforth assume \(X\ne V\), so \(1\le d<N\).

\paragraph{A potential on feasible spanning trees.}
Call a spanning tree \(S\) \emph{feasible} if its complementary edges admit an orientation in which every vertex has the same outdegree as in the input orientation \(D\). We define a potential on feasible spanning trees. Give weight \(d\) to each vertex of \(X\) and weight two to every other vertex. Let \(w(U)\) denote the total weight of a vertex set \(U\). For a feasible spanning tree \(S\) and an edge \(e\in E(S)\), let \(U\) and \(V\setminus U\) be the vertex sets of the two components of \(S-e\). Define the \emph{charge} \(\chi_S(e)\) of \(e\) and the potential \(\phi(S)\) by
\[
\chi_S(e)=\min\{d,w(U),w(V\setminus U)\},
\qquad
\phi(S)=\sum_{e\in E(S)}\chi_S(e).
\]
An edge of \(S\) is \emph{saturated} if its charge is \(d\), and \emph{light} otherwise. Since every charge is at most \(d\), a saturated edge cannot gain charge. Since \(w(V)>2d\), every light edge \(e\) has a unique side \(A_e\) of weight less than \(d\), called its small side. As each vertex in \(X\) has weight \(d\), \(A_e\) is disjoint from \(X\) and \(\chi_S(e)=2|A_e|\).

\paragraph{Exchanges preserving the outdegrees.}
For a feasible spanning tree \(S\), fix a complementary orientation \(D_S\) witnessing its feasibility. For each vertex \(x\in V\), let \(U_x=\{u:(u,x)\in E(D_S)\}\). Root \(S\) at \(x\), and simultaneously replace the parent edge of each \(u\in U_x\) by \(\{u,x\}\). Let \(S_x\) be the resulting graph. It is a tree: following parent edges from any vertex follows its path in $S$ toward \(x\) until it reaches \(x\) or a vertex in \(U_x\), which now points directly to \(x\).

For each \(u\in U_x\), replace the complementary arc \((u,x)\) by the arc from \(u\) to its parent in $S$. In the complementary orientation, each affected vertex loses and gains one outgoing arc, so the resulting tree is feasible. Similarly, a leaf can be detached and reattached to any of its outneighbors in \(D_S\) while preserving every outdegree, and hence feasibility.

Start with \((S,D_S)=(T,D)\). In each round, choose a diametral path \(P\) of \(S\). Test the reattachment of each endpoint of \(P\) to each of its outneighbors in \(D_S\), as well as the exchange \(S_x\) for every \(x\in V\). If a candidate decreases \(\phi\), accept it together with its complementary orientation and begin a new round. If no candidate decreases $\phi$, terminate the process. Since \(\phi\) is a nonnegative integer, the process terminates. Let \(Q\), \(D_Q\), and \(P\) be the tree, orientation, and diametral path in the final round, and define \(U_x\) using \(D_Q\) from now on. We output \(Q\) and prove that \(\diam(Q)=O(|X|+N/d)\). By the termination condition, neither \(Q_x\) nor a reattachment of either endpoint of \(P\) to an outneighbor in \(D_Q\) decreases \(\phi\).

For each \(x\in V\), root both \(Q\) and \(Q_x\) at \(x\). To compare charges, pair each edge of \(Q\) with the edge of \(Q_x\) having the same child. Let \(f_x(e)\) be the paired edge and let
\[
\Delta_x(e)=\chi_{Q_x}(f_x(e))-\chi_Q(e).
\]
The subtree rooted at any vertex in \(Q_x\) is contained in the subtree rooted at the same vertex in \(Q\). For each non-root vertex \(v\), the subtree rooted at \(v\) in \(Q_x\) also excludes the subtree in \(Q\) of every proper descendant of \(v\) that belongs to \(U_x\), whose path to \(x\) in \(Q_x\) bypasses \(v\).

By the stopping condition,
\[
0\le\sum_{x\in V}\bigl(\phi(Q_x)-\phi(Q)\bigr)
=\sum_{e\in E(Q)}\sum_{x\in V}\Delta_x(e).
\]

\paragraph{Bounding the small sides.}
We first prove that the total size of all small sides is at most \(2N\). We will use this estimate to bound the number of small sides containing each endpoint of \(P\). Let
\[
s=\sum_{e\text{ light}}|A_e|.
\]
Fix a light edge \(e\), let \(a=|A_e|\), and let \(v\) be the endpoint of \(e\) in \(A_e\). For \(x\in A_e\), the charge of the paired edge can increase by at most \(d-2a\). For \(x\notin A_e\), the subtree rooted at \(v\) in \(Q\) is exactly \(A_e\). The subtree rooted at \(v\) in \(Q_x\) is contained in \(A_e\) and excludes every vertex in \((A_e\setminus\{v\})\cap U_x\). Therefore, $\Delta_x(e) \le -2|(A_e\setminus\{v\})\cap U_x|$. Each vertex of \(A_e\setminus\{v\}\) has outdegree at least \(d\) in \(D_Q\), so at least \(d-a\) of its outneighbors lie outside \(A_e\). Summing over \(x\in V\) gives
\[
\sum_{x\in V}\Delta_x(e)
\le a(d-2a)-2(a-1)(d-a)
\le-d|A_e|+2d.
\]
Saturated edges have nonpositive changes, and there are fewer than \(N\) light edges. Hence
\[
0\le\sum_{e\in E(Q)}\sum_{x\in V}\Delta_x(e)\le-ds+2dN,
\qquad\text{so}\qquad s\le2N.
\]

\paragraph{Bounding the saturated middle.}
We next show that the saturated edges on \(P\) form a subpath \(H\) with at most \(3|X|+9N/d\) vertices. They are consecutive because, along \(P\), the weight on one side of successive cuts increases while the weight on the other side decreases. Let \(H\) be this subpath; if there is no saturated edge, let \(H\) consist of one vertex of \(P\). Let \(\ell=|V(H)|\). For each \(v\in V(H)\), let \(F_v\) be the vertex set of the component of \(Q-E(H)\) containing \(v\). These sets partition \(V\), and
\[
\sum_{v\in V(H)}w(F_v)=w(V)\le2N+d|X|.
\]
For this counting argument, mark exactly \(d\) outgoing arcs in \(D_Q\) at each vertex outside \(X\), and none at vertices of \(X\), and keep these marks fixed. Fix a vertex \(x\), and let \(z\in V(H)\) be the unique vertex with \(x\in F_z\). Consider a marked arc \((u,x)\) with \(u\in V(H)\) and \(\operatorname{dist}_H(u,z)\ge2\). The first two edges of the \(u\)--\(x\) path lie in \(H\). Root \(Q\) at \(x\), let \(v\) be the parent of \(u\), and let \(e\) be the parent edge of \(v\). The subtree rooted at \(v\) in \(Q\) consists of \(F_v\) and the subtree rooted at \(u\) in \(Q\). Redirecting \(u\) to \(x\) removes the latter subtree, and other redirections can only shrink what remains. Thus the subtree at \(v\) in \(Q_x\) is contained in \(F_v\), even if \(v\) is also redirected. Since \(e\) is saturated,
\[
\Delta_x(e)\le w(F_v)-d.
\]
\begin{figure}[htbp]
\centering
\begin{tikzpicture}[
  vertex/.style={circle,draw=black!65,fill=white,inner sep=0pt,minimum size=6mm,font=\small},
  branch/.style={circle,draw=black!65,fill=white,inner sep=0pt,minimum size=2mm},
  tree/.style={-{Stealth[length=2mm]},draw=black!65,line width=.8pt},
  marked/.style={-{Stealth[length=2.2mm]},draw=blue!65!black,line width=1pt,densely dashed},
  every node/.style={font=\small}
]
\fill[black!6,rounded corners=7pt] (1.45,-1.15) rectangle (3.15,.4);
\fill[green!8,rounded corners=7pt] (6.35,-1.15) rectangle (9.05,.4);
\node[vertex] (u) at (0,0) {\(u\)};
\node[vertex] (v) at (2.3,0) {\(v\)};
\node[vertex] (r) at (4.6,0) {};
\node[vertex] (z) at (6.9,0) {\(z\)};
\node[vertex] (x) at (8.5,-.75) {\(x\)};
\node[branch] (b1) at (1.85,-.8) {};
\node[branch] (b2) at (2.75,-.8) {};
\draw[tree] (u)--(v);
\draw[tree,line width=1.4pt] (v)--(r);
\draw[tree,dashed] (r)--(z);
\draw[tree] (z)--(x);
\draw[tree] (b1)--(v);
\draw[tree] (b2)--(v);
\draw[marked] (u.north) .. controls (2.1,1.9) and (8.5,1.9) .. (x.north);
\node[above,text=blue!65!black] at (4.25,1.35) {marked arc \((u,x)\)};
\node[below] at (1.15,-.3) {in \(H\)};
\node[above] at (3.45,.25) {\(e\), saturated};
\node[below] at (5.75,-.3) {path in \(H\)};
\node[below] at (2.3,-1.15) {\(F_v\)};
\node[below] at (7.7,-1.15) {\(F_z\)};
\end{tikzpicture}
\caption{The shaded sets \(F_v\) and \(F_z\) are the vertex sets of components of \(Q-E(H)\), with \(x\in F_z\). The first two edges of the \(u\)--\(x\) path in \(Q\) lie in \(H\). Redirecting \(u\) to \(x\) removes the subtree rooted at \(u\) in \(Q\) from the subtree rooted at \(v\). The subtree at \(v\) in \(Q_x\), below the edge paired with \(e\), is therefore contained in \(F_v\). Straight arrows indicate the temporary rooting at \(x\), not the complementary orientation.}
\label{fig:legalpair}
\end{figure}
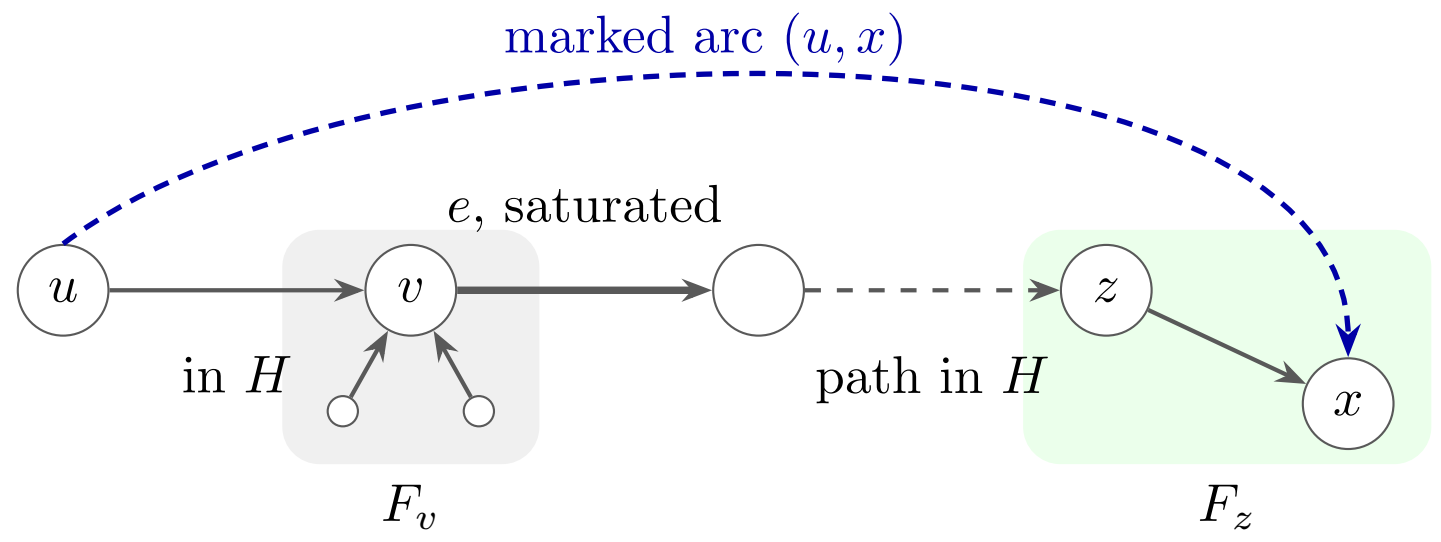

There are at least \(d(\ell-|X|)-3N\) such marked arcs. Indeed, vertices of \(V(H)\setminus X\) supply \(d\) marks each. For a fixed vertex \(x\), at most three vertices of \(H\) lie within distance one of \(z\), and simplicity allows at most one arc from each of them to \(x\). These exclusions account for at most \(3N\) marks in total.

For a fixed \(x\), different marked arcs \((u,x)\) with \(u\in V(H)\) and \(\operatorname{dist}_H(u,z)\ge2\) select different edges \(e\): the edge \(e\), oriented toward \(x\), determines its child \(v\) and the unique preceding neighbor \(u\) on \(H\). Across all \(x\in V\), each \(F_v\) is used at most \(2d\) times: an associated marked arc has its tail at one of the at most two neighbors of \(v\) on \(H\), and each such tail has only \(d\) marks. All other saturated-edge changes are nonpositive. Consequently,
\[
\begin{aligned}
\sum_{e\text{ saturated}}\sum_{x\in V}\Delta_x(e)
&\le-d\bigl(d(\ell-|X|)-3N\bigr)+2d\,w(V)\\
&\le-d^2\ell+3d^2|X|+7dN.
\end{aligned}
\]
Adding the light-edge contribution, which is at most \(2dN\), and using the stopping condition gives
\[
0\le-d^2\ell+3d^2|X|+9dN,
\qquad\text{hence}\qquad
\ell\le3|X|+\frac{9N}{d}.
\]

\paragraph{Bounding the light ends.}
It remains to bound the number of light edges on \(P\). Let \(r(v)\) count the small sides \(A_e\) containing \(v\). Double counting the pairs $(e,v)$ such that $e$ is a light edge and $v\in A_e$ gives \(\sum_{v\in V} r(v)=s\le2N\). We show that each endpoint \(a\) of \(P\) satisfies \(r(a)\le2N/d+1\).

If \(r(a)=0\), there is nothing to prove. Otherwise, $a$ belongs to a small side $A_e$ for some $e$. Since $w(A_e)<d$, we have \(a\notin X\) and \(d>2\). As an endpoint of a diametral path, \(a\) is a leaf. Reattach \(a\) to any of its outneighbors \(x\) in \(D_Q\), obtaining \(\widehat Q\). The leaf edges incident with \(a\) in \(Q\) and \(\widehat Q\) both isolate \(a\), so both have charge two. Among the other light cuts, a small side containing \(a\) but not \(x\) loses two units of weight, while one containing \(x\) but not \(a\) gains at most two. All other light-cut charges are unchanged, and saturated cuts cannot gain. Excluding the leaf cut at \(a\) in \(Q\) from the count therefore gives
\[
0\le\phi(\widehat Q)-\phi(Q)
\le2\bigl(r(x)-r(a)+1\bigr).
\]
The first inequality follows from the stopping condition. Thus every outneighbor \(x\) of \(a\) in \(D_Q\) has \(r(x)\ge r(a)-1\). There are at least \(d\) such vertices, so
\[
d(r(a)-1)\le \sum_{(a,x)\in E(D_Q)} r(x)\le \sum_{v\in V} r(v)=s\le2N.
\]
The same bound holds for the other endpoint \(b\).

Every light edge of \(P\) has exactly one of \(a,b\) in its small side since each cut on $P$ separates $a,b$. There are therefore at most \(r(a)+r(b)\) such edges, whereas the saturated edges are the \(\ell-1\) edges of \(H\). Combining the bounds already proved yields
\[
\diam(Q)\le\ell-1+r(a)+r(b)
\le3|X|+\frac{13N}{d}+1
=O(|X|+N/d).
\]
The orientation \(D_Q\) preserves all input outdegrees by feasibility, so it serves as \(D'\).
\end{proof}

\subsection{Algorithm and running time}
\label{subsec:outdegree-shortening-algorithm}

Let \(m=|E|\). The construction takes \(O((dN+1)(m+N^2))=O(dN^3)\) time in the nontrivial case. The case \(X=V\) returns the input pair after \(O(m+N)\) preprocessing.

\paragraph{Number of rounds.}
Use fixed orders to resolve all choices. The integer potential is nonnegative and at most \(d(N-1)\), so there are at most \(d(N-1)\) accepted exchanges, followed by one final round.

\paragraph{Cost of a round.}
Two tree traversals find a diametral path in \(O(N)\) time. Scanning the current orientation groups all arcs by head in \(O(m+N)\) time. There are at most \(N+2(N-1)=O(N)\) candidates. For a candidate \(S_x\), root \(S\) at \(x\) and replace the specified parent edges; for an endpoint candidate, replace its single leaf edge. Each candidate tree can be constructed in \(O(N)\) time. Another traversal computes all subtree weights and hence every cut charge, also in \(O(N)\) time.

Rejected candidates need no copy of the complementary orientation. For an accepted candidate, update its exchanged edge directions in \(O(m+N)\) time. A round therefore costs \(O(m+N^2)\). Multiplying by the bound on the number of rounds proves the claimed running time, including initialization and the final unsuccessful round.

\section{Proof of the Main Theorem}
\label{sec:combine}

We combine the results of \cref{sec:leaf-supplying,sec:improvement,sec:outdegree-shortening} as in Steps 1--4 of \cref{sec:overview} to obtain $k$ edge-disjoint spanning trees with bounds in terms of $k$ (Proposition~\ref{prop:k-bounds}). We then strengthen these bounds in terms of the minimum degree (Lemma~\ref{cor:minimum-degree}).

\begin{proposition}\label{prop:k-bounds}
Let $G$ be an $n$-vertex simple graph with $n\ge2$. Given $k\ge1$ edge-disjoint spanning trees of $G$, we can construct $k$ edge-disjoint spanning trees, each of diameter $O(n/k)$ and with $O((n/k)\log(k+1))$ internal vertices, in expected
$
\tilde{O}(k^{7/2}n^{3/2}+kn^3)
$
time.
\end{proposition}

The following lemma strengthens the bounds in \cref{prop:k-bounds} by replacing $k$ with the minimum degree $\delta$, while retaining all $k$ trees. We partition $G$ into $O(n/\delta)$ parts, construct trees within each part using \cref{prop:k-bounds}, and connect them using the given edge-disjoint spanning trees.

\begin{lemma}\label{cor:minimum-degree}
Let $G$ be an $n$-vertex simple graph with minimum degree $\delta\ge2$. Given $k\ge1$ edge-disjoint spanning trees of $G$, we can construct $k$ edge-disjoint spanning trees, each of diameter $O(n/\delta)$ and with $O((n/\delta)\log\delta)$ internal vertices, in expected
$
\tilde{O}\!\left(n^4+\delta^{7/2}n^{3/2}+\delta n^3\right)
$
time.
\end{lemma}

\begin{proof}
If $\delta<16k$, the result follows from \cref{prop:k-bounds} when $k\ge2$; when $k=1$, any spanning tree suffices because $2\le\delta<16$. Assume $\delta\ge16k$, and let $t=\lfloor\delta/8\rfloor$. Then $t\ge2k\ge2$ and $t=\Theta(\delta)$.

Initially, let $X=V(G)$. If $X$ can be partitioned into two nonempty sets with fewer than $2t$ edges of $G[X]$ joining them, split $X$ into these sets and apply the same procedure recursively to each. Let $V_1,\ldots,V_q$ be the sets obtained when no further split is possible. Let $M$ count the edges joining distinct parts. There are $q-1$ splits, each contributing at most $2t$ such edges, so
\[
M\le2t(q-1).
\]
For a vertex set $X$, let $\partial X$ be the set of edges with exactly one endpoint in $X$. Any vertex in $V_j$ has at most $|V_j|-1$ neighbors inside $V_j$ and at most $|\partial V_j|$ neighbors outside it. Thus $\delta+1\le |V_j|+|\partial V_j|$ for every $j$. Summing $\delta+1\le |V_j|+|\partial V_j|$ over $j\in[q]$, and noting that each of the $M$ edges joining distinct sets is counted twice, gives
\[
q(\delta+1)\le n+2M\le n+4t(q-1)\le n+\frac{\delta}{2}q.
\]
Therefore, $q\le2n/\delta$.

For every $j\in[q]$ with $|V_j|\ge2$, the graph $G[V_j]$ is $2t$-edge-connected. By the Nash-Williams--Tutte theorem, it contains $t$ edge-disjoint spanning trees. Applying \cref{prop:k-bounds}, we obtain $t$ such trees with the stated diameter and internal-vertex bounds. Retain any $k$ of these trees and denote them by $S_{1,j},\ldots,S_{k,j}$. Since $t=\Theta(\delta)$, for every $i\in[k]$,
\[
\diam(S_{i,j})=O(|V_j|/\delta),
\qquad
|\int(S_{i,j})|=O\!\left(\frac{|V_j|}{\delta}\log\delta\right).
\]
If $V_j$ is a singleton, let each $S_{i,j}$ be its one-vertex tree.

Let $T'_1,\ldots,T'_k$ be the given edge-disjoint spanning trees of $G$. For every $i\in[k]$, contract each vertex set $V_j$ to a single vertex in $T'_i$ and discard loops. The resulting multigraph is connected. Choose a spanning tree of it, and let $F_i$ be the corresponding $q-1$ edges of $T'_i$. Then
\[
T_i=\left(V(G),\ F_i\cup\bigcup_{j=1}^{q}E(S_{i,j})\right)
\]
is a spanning tree of $G$. The sets $F_i$ are pairwise disjoint because the trees $T'_i$ are edge-disjoint. Each edge in $F_i$ joins distinct sets $V_j$, whereas each local tree lies within one such set. Since the local trees within each $V_j$ are also edge-disjoint, so are $T_1,\ldots,T_k$.

For every $i\in[k]$, a path in $T_i$ visits each vertex set $V_j$ at most once because contracting these sets gives a tree. Hence
\[
\diam(T_i)\le\sum_{j=1}^{q}\diam(S_{i,j})+(q-1)=O(n/\delta).
\]
A vertex that becomes internal when the edges of $F_i$ are added must be an endpoint of one of these edges. Since $|F_i|=q-1$, we obtain
\[
|\int(T_i)|
\le\sum_{j:\,|V_j|\ge2}|\int(S_{i,j})|+2(q-1)
=O\!\left(\frac n\delta\log\delta\right).
\]
Thus the trees $T_1,\ldots,T_k$ satisfy both required bounds.

We now bound the running time. If $\delta<16k$, one application of Proposition~\ref{prop:k-bounds} suffices. Otherwise, the recursive partition uses $O(n/\delta)$ tests for a cut of size less than $2t$. To perform each test, fix a vertex of the induced subgraph and compute a flow to every other vertex, stopping each flow computation after finding $2t$ units of flow. If any flow computation terminates with value less than $2t$, it gives a cut of size less than $2t$. If every flow reaches $2t$, no such cut exists. An augmenting-path implementation takes $O(nt(n+m))$ time per induced subgraph, where $m=|E(G)|$. Since $t=O(\delta)$ and $m=O(n^2)$, the partition and construction of the induced subgraphs take $O(n^4)$ time in total.

For each non-singleton $V_j$, find $t$ edge-disjoint spanning trees of $G[V_j]$ in $\tilde{O}(t|V_j|^2+|E(G[V_j])|)$ time, and apply Proposition~\ref{prop:k-bounds} to them. Since $\sum_j|V_j|\le n$, all these applications take expected
\[
\tilde{O}\!\left(t^{7/2}n^{3/2}+tn^3\right)
\]
time in total. The final trees can be assembled in $O(kn)$ time from the local trees and the given trees. As $t=O(\delta)$, the stated bound follows. The expected time includes the cost of rejected samples and retries.
\end{proof}

\begin{proof}[Proof of \cref{thm:main}]
The partition condition is necessary by the Nash-Williams--Tutte theorem. If it holds, the same theorem gives $k$ edge-disjoint spanning trees. When $\delta\ge2$, Lemma~\ref{cor:minimum-degree} constructs trees satisfying both minimum-degree bounds. When $\delta=1$, the given trees satisfy these bounds. Finding the initial trees takes $\tilde{O}(n^3)$ time, and the subsequent construction takes expected $\tilde{O}(n^5)$ time.
\end{proof}

\section{Hardness results}
\label{sec:hardness-results}

This section proves two hardness results for \(n\)-vertex simple graphs
\(G\) promised to contain \(k\) edge-disjoint spanning trees.
\begin{itemize}
\item For every fixed \(k\ge2\), it is NP-hard to decide whether \(G\) contains
\(k\) edge-disjoint spanning trees, each of diameter at most
\(\frac12 \cdot\frac{n}{\delta(G)}\).
\item For every sufficiently large fixed \(k\), it is NP-hard to decide
whether \(G\) contains \(k\) edge-disjoint spanning trees, each with at most
\(\frac{1}{600}\cdot\frac{n \ln\delta(G)}{\delta(G)}\) internal vertices.
\end{itemize}

The following two observations give the edge-disjoint trees used in both reductions.

\begin{observation}
\label[observation]{lem:rooted-clique-packing}
Let \(C\) be the clique on
\(\{r,a_1,b_1,\ldots,a_k,b_k\}\). Then \(C\) contains \(k\)
pairwise edge-disjoint spanning trees \(H_1,\ldots,H_k\), each of
depth at most two from \(r\), such that
\(N_{H_i}(r)=\{a_i,b_i\}\) for every \(i\in[k]\).
\end{observation}

\begin{proof}
For each \(i\), start \(H_i\) with \(ra_i,rb_i\). For each \(i<j\),
put \(a_i a_j,b_i b_j\) in \(H_i\), and \(a_i b_j,b_i a_j\) in
\(H_j\). In \(H_i\), every vertex outside \(\{r,a_i,b_i\}\) is a
leaf adjacent to \(a_i\) or \(b_i\). Thus \(H_i\) has the stated
properties. The root edges are distinct across the trees, and the four
edges between any two designated pairs are used once each. Hence the
trees are pairwise edge-disjoint.
\end{proof}

\begin{observation}
\label[observation]{lem:clique-double-stars}
Let \(C\) be a clique with \(k\) disjoint pairs of vertices
\((x_1,y_1),\ldots,(x_k,y_k)\). Then \(C\) contains \(k\) pairwise
edge-disjoint spanning trees \(F_1,\ldots,F_k\) such that \(x_i y_i\in E(F_i)\)
and every vertex outside \(\{x_i,y_i\}\) is a leaf of \(F_i\).
\end{observation}

\begin{proof}
For each \(i\), start \(F_i\) with \(x_i y_i\). For each \(j>i\), add
\(x_i x_j,y_i y_j\) to \(F_i\) and \(x_i y_j,y_i x_j\) to \(F_j\).
Attach every vertex outside the designated pairs to \(x_i\) in \(F_i\).
Each \(F_i\) is a spanning double star centered at \(x_i,y_i\). The four
edges between any two designated pairs are used once in total, while all
other edges used by the trees are distinct. Thus the trees are pairwise
edge-disjoint.
\end{proof}

\subsection{Hardness of finding spanning trees with small diameter}
\label{sec:diameter-hardness}

We first show that finding edge-disjoint spanning trees of depth two
from a prescribed root is hard, even when the input is promised to
contain the required number of edge-disjoint spanning trees. We then
attach a long tail to turn this rooted depth requirement into a diameter
bound.

\citet{bermond1992broadcasting} proved the
following lemma for \(k=2\), using a reduction due to Alon. We generalize
their reduction to every fixed \(k\ge2\).

\begin{lemma}
\label{lem:rooted-depth-two-hardness}
For every fixed integer \(k\ge2\), it is NP-hard to decide, given a
simple graph \(G\) containing \(k\) pairwise edge-disjoint spanning
trees and a specified root \(r\), whether \(G\) has \(k\) pairwise
edge-disjoint spanning trees of depth at most two from \(r\).
\end{lemma}

\begin{proof}
Fix \(k\ge2\). We reduce from monotone NAE-3-SAT, which is
NP-hard~\citep{lovasz1973coverings}. An instance \(\Phi\) has variable
set \(X\) and clause set \(\mathcal C\), with each clause containing
three distinct unnegated variables. The question is whether each clause
can contain both a true and a false variable.

Create a vertex for each variable, a root \(r\), and vertices
\(a_i,b_i\) for \(i\in[k]\); make these vertices a clique \(Q\). For
each clause \(C\), add a vertex \(z_C\) adjacent exactly to its three
variable vertices and to \(a_3,\ldots,a_k\). Add a vertex \(z_0\)
adjacent exactly to \(a_1,\ldots,a_k\). Let \(G\) be the resulting
simple graph. Its size is polynomial in the size of \(\Phi\).

First, \(G\) satisfies the promise, regardless of the answer for
\(\Phi\). Apply Observation~\ref{lem:rooted-clique-packing} to the
clique on \(\{r\}\cup\{a_i,b_i:i\in[k]\}\). In tree \(i\), attach
each variable and \(z_0\) to \(a_i\). For each clause \(C\), choose
two distinct variables in \(C\), and attach \(z_C\) through the first
in tree 1, the second in tree 2, and \(a_i\) in each tree \(i\ge3\).
This gives \(k\) pairwise edge-disjoint spanning trees of \(G\).

Suppose \(\Phi\) has an NAE assignment; color true variables red and
false variables blue. Start again with the trees from the observation.
For each red variable \(x\), attach \(x\) to \(r\) in tree 1 and to
\(a_i\) in every tree \(i\ne1\). For each blue variable \(x\), attach
\(x\) to \(r\) in tree 2 and to \(a_i\) in every tree \(i\ne2\).
Attach \(z_0\) through \(a_i\) in tree \(i\). Finally, attach each
\(z_C\) through a red variable of \(C\) in tree 1, a blue variable in
tree 2, and \(a_i\) in every tree \(i\ge3\). The resulting trees are
edge-disjoint and have depth at most two from \(r\).

Conversely, suppose \(G\) has \(k\) edge-disjoint spanning trees of
depth at most two from \(r\). Since \(z_0\) is not adjacent to \(r\)
and has precisely the neighbors \(a_1,\ldots,a_k\), its paths from
\(r\) in the trees use these \(k\) neighbors distinctly. Relabel the
trees so that tree \(i\) contains \(ra_i z_0\). For each clause \(C\),
the paths from \(r\) to \(z_C\) in trees 1 and 2 must pass through
variables of \(C\): a path through \(a_j\), \(j\ge3\), would use
\(ra_j\), already used by tree \(j\). Color \(x\) red if \(rx\) belongs
to tree 1, and blue otherwise. The variable on the tree-1 path is red,
whereas the variable on the tree-2 path is blue because the trees are
edge-disjoint. Thus every clause contains both colors, giving an NAE
assignment. This proves the reduction.
\end{proof}

\begin{theorem}
\label{thm:diameter-hardness}
For every fixed integer \(k\ge2\), it is NP-hard to decide, given an
\(n\)-vertex simple graph \(G\) containing \(k\) pairwise edge-disjoint
spanning trees, whether \(G\) has \(k\) pairwise edge-disjoint spanning
trees of diameter at most \(n/(2\delta(G))\).
\end{theorem}

\begin{proof}
Fix \(k\ge2\), and let \((F,r)\) be an instance of
Lemma~\ref{lem:rooted-depth-two-hardness}; in particular, \(F\) contains
\(k\) edge-disjoint spanning trees. Put \(N=|V(F)|\). We attach a tail
at \(r\): its far endpoint forces every vertex of \(F\) to remain close
to \(r\) in any spanning tree of small enough diameter.

For a set \(D\), define a block \(B(D)\) with terminals \(\alpha,\beta\),
vertices \(a_i,b_i\) for \(i\in[k]\), and dummy vertices in \(D\), all
distinct. Its nonterminal vertices form a clique; the terminal edges are
exactly \(\alpha a_i\) and \(\beta b_i\) for \(i\in[k]\).

Let \(\ell\ge1\) and \(D\) be chosen later. Form a tail \(H\) from
one copy of \(B(D)\) followed by \(\ell-1\) copies of
\(B(\emptyset)\), initially disjoint from each other and from \(F\).
Identify consecutive right and left terminals, and let \(s,t\) be the
first and last terminals of the tail. Identify \(s\) with \(r\) to
obtain \(G\); see Figure~\ref{fig:diameter-tail}.

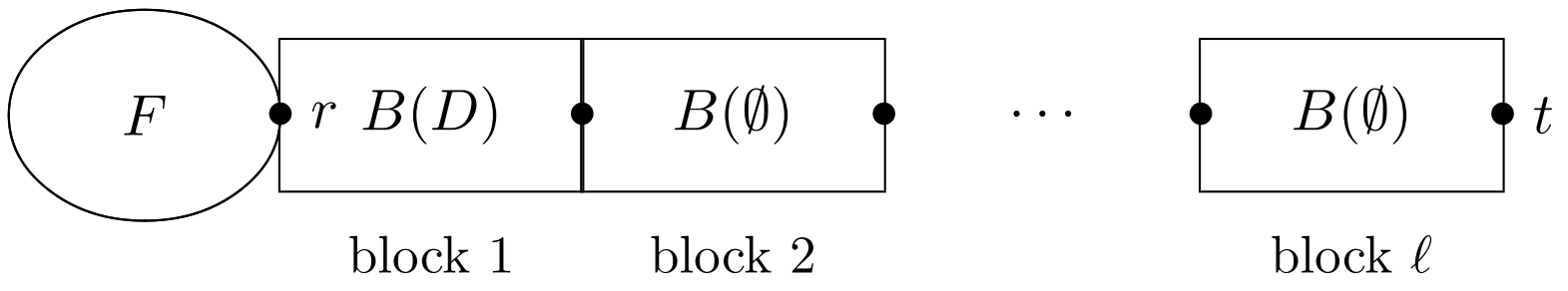
\begin{figure}[htbp]
\centering
\setlength{\unitlength}{1mm}%
\begin{picture}(113,22)
  \qbezier(18,11)(18,13.90)(15.36,15.95)
  \qbezier(15.36,15.95)(12.73,18)(9,18)
  \qbezier(9,18)(5.27,18)(2.64,15.95)
  \qbezier(2.64,15.95)(0,13.90)(0,11)
  \qbezier(0,11)(0,8.10)(2.64,6.05)
  \qbezier(2.64,6.05)(5.27,4)(9,4)
  \qbezier(9,4)(12.73,4)(15.36,6.05)
  \qbezier(15.36,6.05)(18,8.10)(18,11)
  \put(9,11){\makebox(0,0){\(F\)}}
  \put(18,6){\framebox(20,10){}}
  \put(38,6){\framebox(20,10){}}
  \put(79,6){\framebox(20,10){}}
  \put(18,11){\makebox(0,0){\(\bullet\)}}
  \put(38,11){\makebox(0,0){\(\bullet\)}}
  \put(58,11){\makebox(0,0){\(\bullet\)}}
  \put(79,11){\makebox(0,0){\(\bullet\)}}
  \put(99,11){\makebox(0,0){\(\bullet\)}}
  \put(20,11){\makebox(0,0)[l]{\(r\)}}
  \put(28,11){\makebox(0,0){\(B(D)\)}}
  \put(48,11){\makebox(0,0){\(B(\emptyset)\)}}
  \put(68.5,11){\makebox(0,0){\(\cdots\)}}
  \put(89,11){\makebox(0,0){\(B(\emptyset)\)}}
  \put(101,11){\makebox(0,0)[l]{\(t\)}}
  \put(28,3){\makebox(0,0)[t]{\small block 1}}
  \put(48,3){\makebox(0,0)[t]{\small block 2}}
  \put(89,3){\makebox(0,0)[t]{\small block \(\ell\)}}
\end{picture}\par\vspace{2mm}
\caption{The graph \(G\), obtained by identifying the left end \(s\)
of the tail \(H\) with \(r\).}
\label{fig:diameter-tail}
\end{figure}

\Needspace{11\baselineskip}
\begin{numberedclaim}
\label{clm:diameter-tail}
We have \(\operatorname{dist}_H(s,t)=3\ell\). Moreover, \(H\) contains \(k\)
pairwise edge-disjoint spanning trees \(T_1,\ldots,T_k\) with
\(\operatorname{dist}_{T_i}(s,t)=\operatorname{diam}(T_i)=3\ell\) for every \(i\).
\end{numberedclaim}

\begin{claimproof}
The terminal neighborhoods in a block are disjoint, and
\(\alpha a_i b_i\beta\) is a path, so each block has terminal distance
three and \(\operatorname{dist}_H(s,t)=3\ell\).

For each block, apply Observation~\ref{lem:clique-double-stars} to
its internal clique with the designated pairs \((a_i,b_i)\), and add
\(\alpha a_i,\beta b_i\) to tree \(i\).
The block trees are edge-disjoint and each has diameter three, with
terminal distance three. Concatenating the \(i\)-th trees over the blocks
gives a spanning tree \(T_i\) of \(H\). Its \(s\)-\(t\) path has length
\(3\ell\). Every path with endpoints in different blocks traverses at
most three edges in each block it meets, and a path within one block has
length at most three. Thus every path has length at most \(3\ell\), and the trees have the claimed diameters and are edge-disjoint.
\end{claimproof}

By the promise on \(F\) and Claim~\ref{clm:diameter-tail}, gluing
corresponding trees at \(r=s\) gives \(k\) edge-disjoint spanning trees
of \(G\). Thus the constructed graph satisfies the theorem's promise
and \(\delta(G)\ge k\). The vertex \(t\) has degree \(k\), so
\(\delta(G)=k\).

\Needspace{11\baselineskip}
\begin{numberedclaim}
\label{clm:diameter-equivalence}
The graph \(F\) has \(k\) pairwise edge-disjoint spanning trees of
depth at most two from \(r\) if and only if \(G\) has \(k\) pairwise
edge-disjoint spanning trees of diameter at most \(3\ell+2\).
\end{numberedclaim}

\begin{claimproof}
If \(F\) has the rooted trees, glue them to the trees of
Claim~\ref{clm:diameter-tail}. In each resulting tree, distances within
\(F\) are at most four, within \(H\) at most \(3\ell\), and between
\(F\) and \(H\) at most \(3\ell+2\). Since \(\ell\ge1\), its diameter
is at most \(3\ell+2\).

Conversely, let \(T\) be a spanning tree of \(G\) of diameter at most
\(3\ell+2\). As \(F\) and \(H\) meet only at \(r=s\), the restriction
of \(T\) to each is a spanning tree, and the path from any \(x\in V(F)\)
to \(t\) passes through \(r\). By Claim~\ref{clm:diameter-tail},
\(\operatorname{dist}_T(r,t)\ge \operatorname{dist}_H(r,t)=3\ell\). Therefore
\(\operatorname{dist}_T(x,r)=\operatorname{dist}_T(x,t)-\operatorname{dist}_T(r,t)\le2\). Restricting each of the \(k\)
edge-disjoint trees to \(F\) gives the required rooted trees.
\end{claimproof}

Excluding the dummy vertices in $D$, each block contributes $2k+1$ new vertices. We choose the length and the size of $D$ to make \(n=2k(3\ell+2)\), so the target threshold equals the bound in Claim~\ref{clm:diameter-equivalence}. Choose the smallest \(\ell\ge1\) such that
\(4k+(4k-1)\ell\ge N\), and take
\(|D|=4k+(4k-1)\ell-N\). The number of vertices of \(G\) is
\[
    n=N+(2k+1)\ell+|D|=2k(3\ell+2).
\]
Since \(\delta(G)=k\), the desired threshold is
\(n/(2\delta(G))=3\ell+2\). For fixed \(k\), both \(\ell\) and
\(|D|\) are \(O(N+k)\), so the construction is polynomial in the
input size. Claim~\ref{clm:diameter-equivalence} now maps YES instances
to YES instances and NO instances to NO instances, completing the
reduction.
\end{proof}

\subsection{Hardness of finding spanning trees with few internal vertices}
\label{sec:hardness}

A \emph{Set Cover instance} consists of a nonempty finite universe \(U\)
and nonempty sets \(S_1,\ldots,S_m\subseteq U\) whose union is \(U\);
sets may repeat. A \emph{cover} is an index set
\(J\subseteq[m]\) such that \(\bigcup_{j\in J}S_j=U\). The \emph{frequency} of an element is the
number of sets containing it, counted with multiplicity. 

Our reduction uses the following Set Cover gap. The bound \(|U|\le m\)
will control the number of vertices in the graph, while the frequency
condition will ensure a large minimum degree.
\begin{theorem}
\label{thm:set-cover-hardness}
There is a positive integer \(D_0\) such that, for every fixed integer
\(D>D_0\), the following promise problem is NP-hard. The input is a
Set Cover instance with universe \(U\) and sets \(S_1,\ldots,S_m\), where
\(|U|\le m\) and every element has frequency at least \(D\).
Writing \(\tau\) for its minimum cover size, the task is to distinguish
the following two promised cases:
\[
    \mathrm{YES}:\quad \tau\le\frac mD,
    \qquad
    \mathrm{NO}:\quad \tau>\frac1{100}\frac mD\ln D.
\]
\end{theorem}

\begin{theorem}
\label{thm:internal-vertices-hardness}
For every sufficiently large fixed positive integer \(k\), it is NP-hard
to distinguish the following two cases for an \(n\)-vertex simple graph
\(G\) promised to contain \(k\) pairwise edge-disjoint spanning trees:
\begin{itemize}
    \item[\textnormal{\textsc{Yes}:}] \(G\) contains \(k\) pairwise edge-disjoint
    spanning trees, each with fewer than \(3n/\delta(G)\) internal vertices.
    \item[\textnormal{\textsc{No}:}] Every spanning tree of \(G\) has more than
    \(n\ln\delta(G)/(600\delta(G))\) internal vertices.
\end{itemize}
\end{theorem}

\begin{proof}[Proof of Theorem~\ref{thm:internal-vertices-hardness}]
Fix \(k\) sufficiently large and put \(D=2k\), so that
Theorem~\ref{thm:set-cover-hardness} applies and \(\ln k/600>3\).
The two cases are disjoint for such \(k\), since every graph under the
promise has minimum degree at least \(k\).
Given an instance of that theorem with universe \(U\) and sets
\(S_1,\ldots,S_m\), let \(\tau\) be its minimum cover size and put \(b=m/D\).
We construct an \(n\)-vertex graph \(G\) whose minimum degree is \(D^2\),
whose ratio \(n/\delta(G)\) is comparable to \(b\), and whose spanning
trees encode covers.

For each \(j\in[m]\), create \(D\) vertices \(v_{j,1},\ldots,v_{j,D}\).
Also create vertices \(x_1,y_1,\ldots,x_k,y_k\), and make all these
vertices a clique \(Q\). Add each element \(u\in U\) as a vertex and join
it to every \(v_{j,t}\) such that \(u\in S_j\). These are the only edges incident with
\(U\).

Add a clique \(R\) on new vertices \(a_1,\ldots,a_{D^2+1}\), with exactly
the edges \(a_i x_i\), \(i\in[k]\), between \(R\) and the rest of the
graph. The clique \(R\) fixes the minimum degree at \(D^2\).
Figure~\ref{fig:hardness-construction} illustrates the construction.

\begin{figure}[htbp]
\centering
\setlength{\unitlength}{1mm}%
\begin{picture}(116,59)
  \qbezier(20,36)(20,43.04)(17.07,48.02)
  \qbezier(17.07,48.02)(14.14,53)(10,53)
  \qbezier(10,53)(5.86,53)(2.93,48.02)
  \qbezier(2.93,48.02)(0,43.04)(0,36)
  \qbezier(0,36)(0,28.96)(2.93,23.98)
  \qbezier(2.93,23.98)(5.86,19)(10,19)
  \qbezier(10,19)(14.14,19)(17.07,23.98)
  \qbezier(17.07,23.98)(20,28.96)(20,36)
  \put(10,55){\makebox(0,0){\(R\)}}
  \put(32,20){\framebox(83,32){}}
  \put(73.5,55){\makebox(0,0){\(Q\)}}
  \put(10,45){\line(1,0){25}}
  \put(10,37){\line(1,0){25}}
  \put(10,25){\line(1,0){25}}
  \qbezier(51,31)(53,15.5)(55,0)
  \qbezier(58,31)(56.5,15.5)(55,0)
  \qbezier(71,31)(63,15.5)(55,0)
  \qbezier(51,31)(75.5,15.5)(100,0)
  \qbezier(58,31)(79,15.5)(100,0)
  \qbezier(71,31)(85.5,15.5)(100,0)
  \qbezier(87,31)(93.5,15.5)(100,0)
  \qbezier(94,31)(97,15.5)(100,0)
  \qbezier(108,31)(104,15.5)(100,0)
  \put(10,45){\makebox(0,0){\(\bullet\)}} \put(10,37){\makebox(0,0){\(\bullet\)}}
  \put(10,25){\makebox(0,0){\(\bullet\)}}
  \put(35,45){\makebox(0,0){\(\bullet\)}} \put(42,45){\makebox(0,0){\(\bullet\)}}
  \put(35,37){\makebox(0,0){\(\bullet\)}} \put(42,37){\makebox(0,0){\(\bullet\)}}
  \put(35,25){\makebox(0,0){\(\bullet\)}} \put(42,25){\makebox(0,0){\(\bullet\)}}
  \put(51,31){\makebox(0,0){\(\bullet\)}} \put(58,31){\makebox(0,0){\(\bullet\)}}
  \put(71,31){\makebox(0,0){\(\bullet\)}} \put(87,31){\makebox(0,0){\(\bullet\)}}
  \put(94,31){\makebox(0,0){\(\bullet\)}} \put(108,31){\makebox(0,0){\(\bullet\)}}
  \put(55,0){\makebox(0,0){\(\bullet\)}} \put(100,0){\makebox(0,0){\(\bullet\)}}
  \put(10,47){\makebox(0,0)[b]{\(a_1\)}}
  \put(10,39){\makebox(0,0)[b]{\(a_2\)}}
  \put(26,31){\makebox(0,0){\(\vdots\)}}
  \put(10,27){\makebox(0,0)[b]{\(a_k\)}}
  \put(35,47){\makebox(0,0)[b]{\(x_1\)}}
  \put(42,47){\makebox(0,0)[b]{\(y_1\)}}
  \put(35,39){\makebox(0,0)[b]{\(x_2\)}}
  \put(42,39){\makebox(0,0)[b]{\(y_2\)}}
  \put(35,27){\makebox(0,0)[b]{\(x_k\)}}
  \put(42,27){\makebox(0,0)[b]{\(y_k\)}}
  \put(51,33){\makebox(0,0)[b]{\(\scriptstyle v_{j,1}\)}}
  \put(58,33){\makebox(0,0)[b]{\(\scriptstyle v_{j,2}\)}}
  \put(64.5,31){\makebox(0,0){\(\cdots\)}}
  \put(71,33){\makebox(0,0)[b]{\(\scriptstyle v_{j,2k}\)}}
  \put(87,33){\makebox(0,0)[b]{\(\scriptstyle v_{j',1}\)}}
  \put(94,33){\makebox(0,0)[b]{\(\scriptstyle v_{j',2}\)}}
  \put(101,31){\makebox(0,0){\(\cdots\)}}
  \put(108,33){\makebox(0,0)[b]{\(\scriptstyle v_{j',2k}\)}}
  \put(40,0){\makebox(0,0){\(U\)}}
  \put(55,-2){\makebox(0,0)[t]{\(u_1\)}}
  \put(100,-2){\makebox(0,0)[t]{\(u_2\)}}
\end{picture}\par\vspace{5mm}
\caption{Schematic of \(G\), showing selected vertices and omitting the
internal edges of the cliques \(R\) and \(Q\). The displayed incidences are
\(u_1\in S_j\setminus S_{j'}\) and \(u_2\in S_j\cap S_{j'}\).}
\label{fig:hardness-construction}
\end{figure}
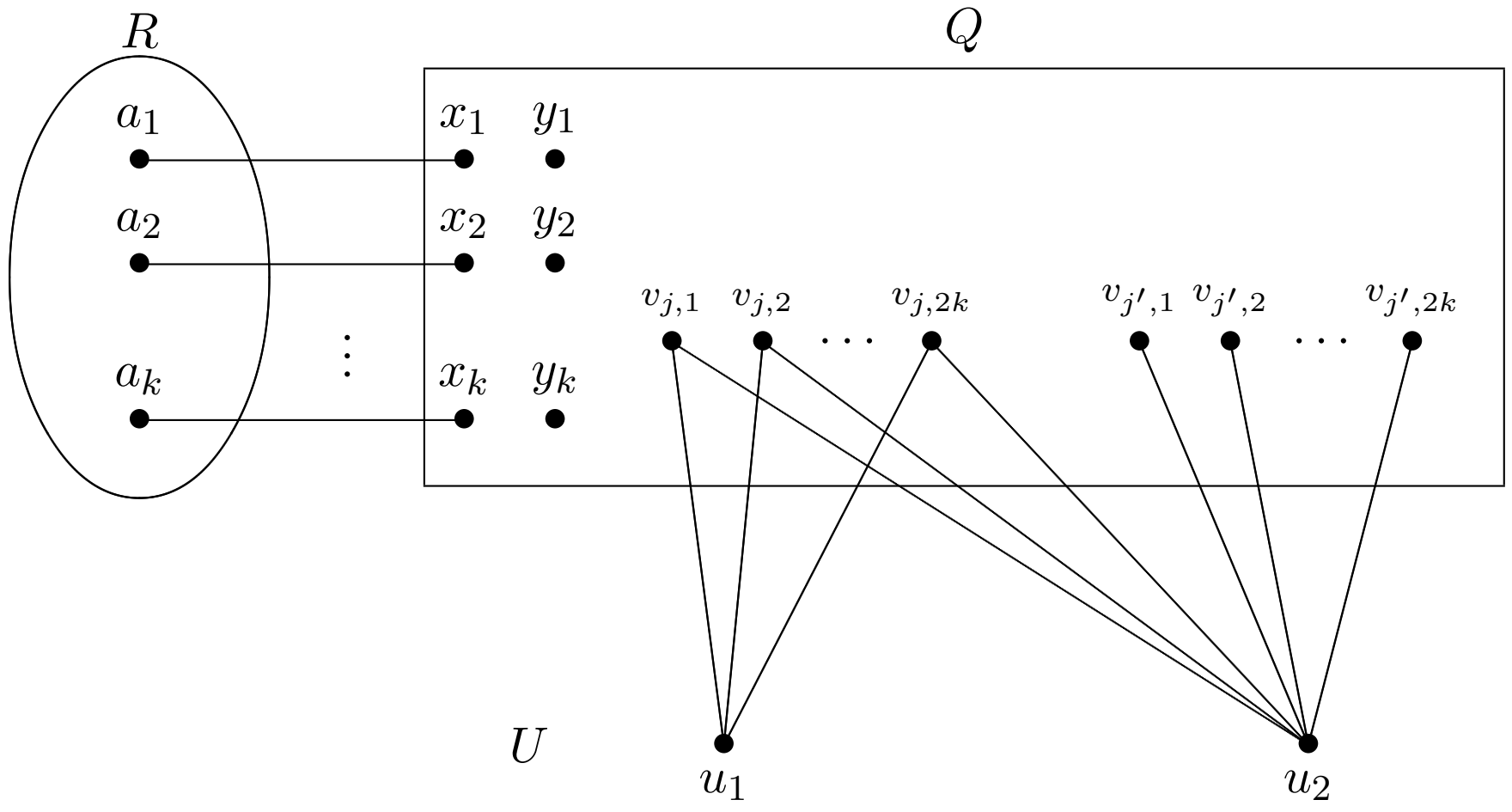

\begin{numberedclaim}
\label{clm:hardness-size-degree}
We have \(b+1<n/\delta(G)\le3b\).
\end{numberedclaim}

\begin{claimproof}
Every vertex \(u\in U\) has degree at least \(D^2\), since its frequency
is at least \(D\).
Since \(U\ne\emptyset\), we have \(m\ge D\), and every vertex of \(Q\)
has degree at least \(Dm+D-1\ge D^2\). In \(R\), the vertices \(a_i\)
with \(i\le k\) have degree \(D^2+1\), while the others have degree \(D^2\).
Therefore \(\delta(G)=D^2\). As \(n=Dm+D+D^2+1+|U|\),
\[
    \frac n{\delta(G)}=b+1+\frac{D+|U|+1}{D^2}>b+1.
\]
Using \(|U|\le m\), \(m\ge D\), and \(k\ge2\), we have
\(D+|U|+1\le Dm\). Hence \(n/\delta(G)\le2b+1\le3b\),
since \(b\ge1\).
\end{claimproof}

\begin{numberedclaim}
\label{clm:hardness-cover-to-trees}
For every cover \(J\), the graph \(G\) has \(k\)
pairwise edge-disjoint spanning trees, each with at most \(|J|+4\)
internal vertices.
\end{numberedclaim}

\begin{claimproof}
Apply Observation~\ref{lem:clique-double-stars} to \(Q\) with pairs
\((x_i,y_i)\), obtaining trees \(A_1,\ldots,A_k\), and to \(R\) with pairs
\((a_i,a_{k+i})\), obtaining trees \(B_1,\ldots,B_k\).
For each \(i\in[k]\), start \(T_i\) with \(A_i\cup B_i\) and the edge
\(x_i a_i\). For each \(u\in U\), choose \(j\in J\) with \(u\in S_j\)
and add the edge \(u v_{j,i}\) to \(T_i\) for every \(i\in[k]\).
Each \(T_i\) is a spanning tree:
the first added edge joins the two clique trees and every other added edge
attaches a new leaf. The trees are pairwise edge-disjoint because the
clique trees are, and their added edges have distinct indices \(i\).
Only \(x_i,y_i,a_i,a_{k+i}\) and vertices \(v_{j,i}\) with \(j\in J\) can
be internal in \(T_i\), so \(|\int(T_i)|\le|J|+4\).
\end{claimproof}

Taking the cover \(J=[m]\) in Claim~\ref{clm:hardness-cover-to-trees}
shows that every graph
produced by the reduction contains \(k\) pairwise edge-disjoint spanning
trees, including in the NO case.

\begin{numberedclaim}
\label{clm:hardness-tree-to-cover}
Every spanning tree \(T\) of \(G\) satisfies \(|\int(T)|\ge\tau\).
\end{numberedclaim}

\begin{claimproof}
Root \(T\) at \(a_1\). The parent of each \(u\in U\) is a vertex \(v_{j,t}\)
with \(u\in S_j\). That parent is internal: it has both a parent of its own
and the child \(u\). The corresponding set indices form a cover \(J\) with
\(|J|\le|\int(T)|\). Hence \(\tau\le|J|\le|\int(T)|\).
\end{claimproof}

In the YES case, \(\tau\le b\), so
Claims~\ref{clm:hardness-size-degree} and~\ref{clm:hardness-cover-to-trees}
give \(k\) pairwise edge-disjoint spanning trees \(T_1,\ldots,T_k\) with
\[
    |\int(T_i)|\le\tau+4\le b+4\le3(b+1)
    <3\frac n{\delta(G)}\qquad(i\in[k]).
\]
In the NO case, \(\tau>\frac b{100}\ln D\), so
Claims~\ref{clm:hardness-size-degree} and~\ref{clm:hardness-tree-to-cover}
give every spanning tree \(T\) the bound
\[
    |\int(T)|\ge\tau
    >\frac b{100}\ln D
    \ge\frac1{300}\frac n{\delta(G)}\ln D
    =\frac1{600}\frac n{\delta(G)}\ln\delta(G).
\]
The construction has polynomial size and takes polynomial time for each
fixed \(k\). This proves the promised hardness gap.
\end{proof}

\section*{AI Methodology}
We used OpenAI’s GPT-6 Astra to investigate whether the result in an earlier version of our manuscript, which gave $(1-\varepsilon)k$ trees satisfying both bounds, could be improved. After extensive exploration, GPT-6 Astra suggested that all $k$ trees could be obtained and provided an approach. We examined this approach and subsequently developed a proof with a substantially different structure. We also used ChatGPT to help check the correctness of our proofs. We are responsible for the final results and their proofs.

\clearpage
\bibliographystyle{cleanrefs}
\bibliography{ref}
\clearpage

\appendix

\section{Proof of the Set Cover hardness theorem}
\label{app:set-cover-hardness}

The following Set Cover hardness result follows from the analysis of
\citet[Section~2.1]{chlebik2008dominating} of the construction of
\citet{trevisan2001bounded}.

\begin{theorem}
\label{thm:bounded-set-cover-source}
There is a positive integer \(B_0\) such that, for every fixed integer
\(B>B_0\), the following promise problem is NP-hard. The input is a
Set Cover instance with universe \(U'\) and sets \(S'_1,\ldots,S'_{m'}\).
An integer \(\beta\) satisfying
\(0\le\beta\le m'\) is also given.\footnote{Chleb\'ik and
  Chleb\'ikov\'a remark that, for every sufficiently large \(B\), it is NP-hard,
  on Set Cover instances satisfying \(\lvert S'_j\rvert\le B\) and
  \(\lvert\{j:u\in S'_j\}\rvert\le(\log B)^{O(1)}\), to distinguish
  \[
      \tau'<t
      \quad\text{from}\quad
      \tau'>t(\ln B-C\ln\ln B),
  \]
  where \(C>0\) is an absolute constant and the threshold \(t\) is part of the
  input, and the sets cover the universe. Choose \(B_0\) so that the gap
  factor exceeds one and \(\ln B/10>1\). Instances with \(t\le0\) are
  trivially NO under the promise. Map them to the one-element, one-set
  instance with \(\beta=0\). For \(t>0\), an empty universe is trivially
  YES; map it to the same instance with \(\beta=1\). Otherwise, delete empty sets; the universe and all
  remaining sets are nonempty, and the optimum is unchanged. Write \(m'\)
  for the number of remaining sets and set
  \(\beta=\min\{\lfloor t\rfloor,m'\}\). Since \(\tau'\) is an
  integer and \(\tau'\le m'\), the YES bound is preserved, while \(\beta\le t\)
  preserves the NO bound. For sufficiently large \(B\), the latter is stronger
  than \(\tau'>(\beta/10)\ln B\), and \((\log B)^{O(1)}\le B\).}
Each set has size at most \(B\), and every element has frequency at most
\(B\). Writing \(\tau'\) for the instance's minimum cover size, the task is to
distinguish the following two promised cases:
\[
    \mathrm{YES}:\quad \tau'\le\beta,
    \qquad
    \mathrm{NO}:\quad \tau'>\frac{\beta}{10}\ln B.
\]
\end{theorem}

\begingroup
\value{theorem}=\getrefnumber{thm:set-cover-hardness}\relax
\setcounter{numberedclaim}{0}
\begin{proof}[Proof of Theorem~\ref{thm:set-cover-hardness}]
Fix an integer \(D>D_0\) and set \(B=\lfloor\sqrt{D/2}\rfloor\).
For sufficiently large \(D_0\), we have \(B>B_0\), so we reduce from the
promise problem in Theorem~\ref{thm:bounded-set-cover-source}.
For an input \(U'\), \(S'_1,\ldots,S'_{m'}\), and \(\beta\), consider
the following linear program (LP) in the variables
\(x_1,\ldots,x_{m'}\):
\begin{equation}
\begin{alignedat}{2}
    0\le x_j &\le 1
        &\qquad &(j=1,\ldots,m'),\\
    \sum_{j:u\in S'_j}x_j &\ge 1
        &&(u\in U'),\\
    \sum_{j=1}^{m'}x_j &\le \beta.&&
\end{alignedat}
\label{eq:fractional-set-cover}
\end{equation}
We test feasibility in polynomial time and, if feasible, obtain an exact
rational solution of polynomial encoding length
\citep{khachiyan1979polynomial}.
This test guides the reduction. If the LP is feasible, its solution tells us
how many copies of each set to use so that every element occurs at least
\(D\) times. If it is infeasible, the source instance must be NO: a cover of
size at most \(\beta\) would give a feasible \(0\)-\(1\) solution. In that
case we output a fixed NO instance instead.

\paragraph{Case (i): the LP is feasible.}
Fix a feasible solution and set \(m=D\beta+2m'\). We will take
\(r_j=\lceil Dx_j\rceil+1\) copies of each \(S'_j\). The following bounds
show that there are at most \(m\) copies and every element belongs to at
least \(D\) of them.

\begin{numberedclaim}
\label{clm:set-cover-multiplicities}
We have
\[
    \sum_{j=1}^{m'}r_j\le m,
    \qquad
    \sum_{j:u\in S'_j}r_j\ge D\quad(u\in U').
\]
\end{numberedclaim}

\begin{claimproof}
Since \(\lceil Dx_j\rceil+1\le Dx_j+2\),
\[
    \sum_{j=1}^{m'}r_j
    \le D\sum_{j=1}^{m'}x_j+2m'
    \le D\beta+2m'=m.
\]
For every \(u\in U'\), the covering constraint gives
\[
    \sum_{j:u\in S'_j}r_j
    \ge D\sum_{j:u\in S'_j}x_j
    \ge D.
\]
\end{claimproof}

Set \(U=U'\). Take \(r_j\) copies of each \(S'_j\), and add copies of
\(S'_1\) until there are exactly \(m\) sets. Denote the resulting sets by
\(S_1,\ldots,S_m\).
Claim~\ref{clm:set-cover-multiplicities} ensures that the extra copies
can be added and that every element has frequency at least \(D\).
All output sets are nonempty. Since \(r_j\ge1\), every source set is
retained, and taking copies does not change the minimum cover size.
Writing \(\tau'\) and \(\tau\) for the source and output minimum cover
sizes, respectively, we have \(\tau=\tau'\).

The feasible LP also bounds the size of the source instance in terms of
\(\beta\). These bounds will verify the output conditions and preserve the
NO gap.

\begin{numberedclaim}
\label{clm:set-cover-source-size}
The source instance satisfies \(|U'|\le B\beta\) and \(m'\le B^2\beta\).
\end{numberedclaim}

\begin{claimproof}
Summing the covering constraints over \(U'\) and using \(|S'_j|\le B\) gives
\[
    |U'|
    \le\sum_{j=1}^{m'}|S'_j|x_j
    \le B\sum_{j=1}^{m'}x_j
    \le B\beta.
\]
Since each \(S'_j\) is nonempty and every element has frequency at most
\(B\), we also have
\[
    m'\le\sum_{j=1}^{m'}|S'_j|
    =\sum_{u\in U'}|\{j:u\in S'_j\}|
    \le B|U'|\le B^2\beta.
\]
\end{claimproof}

Since \(B\le D\), Claim~\ref{clm:set-cover-source-size} gives
\(|U|=|U'|\le B\beta\le D\beta\le m\). Together with
Claim~\ref{clm:set-cover-multiplicities} and nonemptiness of the sets, this
verifies the input conditions of Theorem~\ref{thm:set-cover-hardness}.

If the source instance is YES, then \(\tau'\le\beta\). Since
\(m=D\beta+2m'\), we have \(\beta\le m/D\), and hence
\(\tau=\tau'\le m/D\).

If the source instance is NO, Claim~\ref{clm:set-cover-source-size} and
\(B^2\le D/2\) give
\[
    \frac mD=\beta+\frac{2m'}{D}
    \le \beta+\frac{2B^2\beta}{D}\le2\beta.
\]
For sufficiently large \(D_0\), we have \(\ln B\ge\frac13\ln D\), so
\[
    \tau=\tau'>\frac{\beta}{10}\ln B
    \ge\frac1{20}\frac mD\ln B
    \ge\frac1{60}\frac mD\ln D
    >\frac1{100}\frac mD\ln D.
\]

\paragraph{Case (ii): the LP is infeasible.}
As observed above, the source instance is NO. We output the following fixed
NO instance of Theorem~\ref{thm:set-cover-hardness}.

Let \(q\) be the largest even integer at most \(\log_2D\). Taking
\(D_0\) sufficiently large ensures \(q\ge2\). Define
\[
    U=\{X\subseteq[q]:|X|=q/2\},
    \qquad
    A_j=\{X\in U:j\in X\}\quad(j\in[q]).
\]
Set \(r=\lceil2D/q\rceil\), and take \(r\) copies of each \(A_j\).
Thus \(m=qr\); denote the copies by \(S_1,\ldots,S_m\), and let \(\tau\) be
their minimum cover size.

\begin{numberedclaim}
\label{clm:set-cover-fixed-no}
Every \(S_i\) is nonempty, \(|U|\le m\), and every element has frequency
at least \(D\). Moreover,
\[
    \tau>\frac1{100}\frac mD\ln D.
\]
\end{numberedclaim}

\begin{claimproof}
Every \(A_j\), and hence every \(S_i\), is nonempty. Each \(X\in U\)
belongs to exactly \(q/2\) of the sets \(A_j\),
so it belongs to \(qr/2\ge D\) of their copies. Moreover,
\[
    |U|=\binom q{q/2}\le2^q\le D\le m,
\]
since \(m=qr\ge2D\).

For \(J\subseteq[q]\), the sets \(A_j\) with \(j\in J\) cover \(U\)
exactly when \(|J|>q/2\). Indeed, an uncovered element is a
\(q/2\)-subset of \([q]\setminus J\), and such a subset exists exactly
when \(|J|\le q/2\). Taking copies does not change the minimum cover size, so
\(\tau=q/2+1\). As \(q>\log_2D-2\),
\[
    \tau>\frac12\log_2D>\frac12\ln D.
\]
Finally, \(q\le\log_2D\le D\), so
\(m=q\lceil2D/q\rceil\le2D+q\le3D\). Consequently,
\(\tau>\frac1{100}(m/D)\ln D\).
\end{claimproof}

For sufficiently large \(D_0\), the factor \(\ln D/100\) exceeds one, so
the two target cases are disjoint. The reduction maps YES instances to
YES instances and NO instances to NO instances. It runs in polynomial
time for every fixed \(D>D_0\): in
Case (i), \(m=D\beta+2m'\le(D+2)m'\), while the output in Case (ii) depends
only on \(D\).
\end{proof}
\endgroup

\end{document}